\documentclass[a4paper,11pt,leqno]{amsart}

\usepackage{amsfonts}
\usepackage{amssymb}
\usepackage{amscd}
\usepackage{palatino}
\usepackage{amsthm}
\usepackage{a4wide}

\usepackage[centertags]{amsmath}
\usepackage{eucal,dsfont,accents,bbm,amsfonts,amssymb,amsthm,amsopn}
\usepackage[margin=1.94cm]{geometry}
\usepackage[usenames]{color}
\definecolor{citegreen}{rgb}{0,0.6,0}
\definecolor{refred}{rgb}{0.8,0,0}
\usepackage[colorlinks, citecolor=citegreen, linkcolor=refred]{hyperref}

\theoremstyle{plain}

\newtheorem{teo}{Theorem}[section]

\newtheorem{prop}[teo]{Proposition}
\newtheorem{cor}[teo]{Corollary}
\newtheorem{ackn}{Acknowledgments\!}

\theoremstyle{definition}

\theoremstyle{remark}

\newtheorem{rem}[teo]{Remark}

\newtheorem{note}{Note\!\!\,}

\numberwithin{equation}{section}

\def\SS{{{\mathbb S}}}

\def\RR{{\mathbb R}}

\def\RRR{{\mathrm R}}
\def\WWW{{\mathrm W}}

\def\BBB{{\mathrm B}}

\def\TTT{{\mathrm{T}}}
\def\Ric{{\mathrm {Ric}}}

\def\Scal{{\RRR}}
\def\SSS{{\mathrm S}}
\def\CCC{{\mathrm C}}

\title[The Cotton Tensor and the Ricci Flow]{The Cotton Tensor and the Ricci Flow}

\date{\today}

\author[Carlo Mantegazza]{Carlo Mantegazza}
\address[Carlo Mantegazza]{Scuola Normale Superiore, Piazza
  Cavalieri 7, Pisa, Italy, 56126}
\email[C. Mantegazza]{c.mantegazza@sns.it}

\author[Samuele Mongodi]{Samuele Mongodi}
\address[Samuele Mongodi]{Scuola Normale Superiore, Piazza
  Cavalieri 7, Pisa, Italy, 56126}
\email[S. Mongodi]{samuele.mongodi@gmail.com}

\author[Michele Rimoldi]{Michele Rimoldi}
\address[Michele Rimoldi]{Dipartimento di Matematica e Applicazioni, Universit\'a
  degli Studi di Milano--Bicocca, Via Cozzi 55, Milano, Italy, 20125}
\email[M. Rimoldi]{michele.rimoldi@gmail.com}

\date{\today}

\begin{document}

\begin{abstract}
We compute the evolution equation of the Cotton and the Bach tensor under the Ricci
flow of a Riemannian manifold, with particular attention to the three
dimensional case, and we discuss some applications.
\end{abstract}

\maketitle
\tableofcontents

\section{Preliminaries}

The Riemann curvature operator of a Riemannian manifold $(M^n,g)$ is defined, 
as in~\cite{gahula}, by
$$
\mathrm{Riem}(X,Y)Z=\nabla_{Y}\nabla_{X}Z-\nabla_{X}\nabla_{Y}Z+\nabla_{[X,Y]}Z\,.
$$ 
In a local coordinate system the components of the $(3,1)$--Riemann 
curvature tensor are given by
$\RRR^{l}_{ijk}\frac{\partial}{\partial
  x^{l}}=\mathrm{Riem}\big(\frac{\partial}{\partial
  x^{i}},\frac{\partial}{\partial
  x^{j}}\big)\frac{\partial}{\partial x^{k}}$ and we denote by
$\RRR_{ijkl}=g_{lm}\RRR^{m}_{ijk}$ its $(4,0)$--version.

With the previous choice, for the sphere $\SS^n$ we have
${\mathrm{Riem}}(v,w,v,w)=\RRR_{abcd}v^iw^jv^kw^l>0$. 

\medskip

{\em In all the paper the Einstein convention of summing 
over the repeated  indices will be adopted.}

\medskip

The Ricci tensor is obtained by the contraction 
$\RRR_{ik}=g^{jl}\RRR_{ijkl}$ and $\RRR=g^{ik}\RRR_{ik}$ will 
denote the scalar curvature. 

We recall the interchange of derivative formula,
$$
\nabla^2_{ij}\omega_k-\nabla^2_{ji}\omega_k=\RRR_{ijkp}g^{pq}\omega_q\,,
$$
and Schur lemma, which follows by the second Bianchi identity,
$$
2g^{pq}\nabla_p\RRR_{qi}=\nabla_i\RRR\,.
$$
They both will be used extensively in the computations that follows.

\medskip

The so called Weyl tensor is then 
defined by the following decomposition formula (see~\cite[Chapter~3,
Section~K]{gahula}) in dimension $n\geq 3$,
\begin{equation}\label{decriem}
\RRR_{ijkl}=\frac{1}{n-2}(\RRR_{ik}g_{jl}-\RRR_{il}g_{jk}
+\RRR_{jl}g_{ik}-\RRR_{jk}g_{il})-\frac{\RRR}{(n-1)(n-2)}(g_{ik}g_{jl}-g_{il}g_{jk})
+\WWW_{ijkl}\,.
\end{equation}
The Weyl tensor satisfies all the symmetries of the curvature tensor, moreover, 
all its traces with the metric are zero, 
as it can be easily seen by the above formula.\\
In dimension three $\WWW$ is identically zero for every Riemannian
manifold. It becomes relevant instead when $n\geq 4$ since its
vanishing is a condition equivalent for $(M^n,g)$ to be {\em locally
  conformally flat}, that is, around every point $p\in M^n$ there is a
conformal deformation $\widetilde{g}_{ij}=e^fg_{ij}$ of the original
metric $g$, such that the new metric is flat, 
namely, the Riemann tensor associated to $\widetilde{g}$ is zero in
$U_p$ (here $f:U_p\to\RR$ is a smooth function defined in a open
neighborhood $U_p$ of $p$).\\
In dimension $n=3$, instead, locally conformally flatness is
equivalent to the vanishing of the following {\em Cotton tensor}
\begin{equation}\label{Cottonn}
\CCC_{ijk} =\nabla_k\RRR_{ij} - \nabla_j\RRR_{ik} -
\frac{1}{2(n-1)} \big(\nabla_k\RRR g_{ij} - \nabla_j\RRR g_{ik} \big)\,,
\end{equation}
which expresses the fact that the {\em Schouten} tensor
$$
\SSS_{ij} =\RRR_{ij} - \frac{\RRR g_{ij}}{2(n-1)}
$$ 
is a {\em Codazzi} tensor (see~\cite[Chapter~16, Section~C]{besse}), that is, a symmetric bilinear form $\TTT_{ij}$ such that $\nabla_k\TTT_{ij}=\nabla_i\TTT_{kj}$.

By means of the second Bianchi identity, one can easily get 
(see~\cite{besse}) that
\begin{equation}\label{CottonWeyl}
\,\nabla^l\WWW_{lijk}=-\frac{n-3}{n-2}\CCC_{ijk}\,.
\end{equation}
Hence, when $n\geq 4$, if we assume that the manifold is locally
conformally flat (that is, $\WWW=0$), the Cotton tensor is identically
zero also in this case, but this is only a necessary condition.

By direct computation, we can see that the tensor $\CCC_{ijk}$
satisfies the following symmetries
\begin{equation}\label{CottonSym}
\CCC_{ijk}=-\CCC_{ikj},\,\quad\quad\CCC_{ijk}+\CCC_{jki}+\CCC_{kij}=0\,,
\end{equation}
moreover it is trace--free in any two indices, 
\begin{equation}\label{CottonTraces}
g^{ij}\CCC_{ijk}=g^{ik}\CCC_{ijk}=g^{jk}\CCC_{ijk}=0\,,
\end{equation}
by its skew--symmetry and Schur lemma.

We suppose now that $(M^n, g(t))$ is a Ricci flow in some time
interval, that is, the time--dependent metric $g(t)$ satisfies
\[
\ \frac{\partial}{\partial t}g_{ij}=-2\RRR_{ij}.
\]
We have then the following evolution equations for the Christoffel
symbols, the Ricci tensor and the scalar curvature (see for instance~\cite{hamilton1}),
\begin{eqnarray}\label{evolutioncurv}
\frac{\partial\,}{\partial
  t}\Gamma_{ij}^k&=&-g^{ks}\nabla_i\RRR_{js}-g^{ks}\nabla_j\RRR_{is}+g^{ks}\nabla_s\RRR_{ij}\nonumber\\
\frac{\partial\,}{\partial
  t}\RRR_{ij}&=&\Delta\RRR_{ij}-2\RRR^{kl}\RRR_{kijl}-2g^{pq}\RRR_{ip}\RRR_{jq}\\
\frac{\partial\,}{\partial t}\RRR&=&\Delta\RRR+2|\Ric|^2\,.\nonumber
\end{eqnarray}

\medskip
 
{\em All the computations which follow will be done in a fixed local
  frame, not in a moving frame.}

\bigskip

\begin{ackn} The first and second authors are partially supported by the
Italian FIRB Ideas ``Analysis and Beyond''.
\end{ackn}

\begin{note}
We remark that Huai-Dong Cao also, independently by us, worked out the
computation of the evolution of the Cotton tensor in dimension three,
in an unpublished note.
\end{note}

\bigskip

\section{The Evolution Equation of the Cotton Tensor in 3D}

The goal of this section is to compute the evolution equation under
the Ricci flow of the Cotton tensor $\CCC_{ijk}$ in dimension
three (see~\cite{mancat1} for the evolution of the Weyl tensor), the
general computation in any dimension is postponed to
section~\ref{nsec}.
 
In the special three--dimensional case we have,
\begin{align}
\RRR_{ijkl}=&\,\RRR_{ik}g_{jl}-\RRR_{il}g_{jk}
+\RRR_{jl}g_{ik}-\RRR_{jk}g_{il}-\frac{\RRR}{2}(g_{ik}g_{jl}-g_{il}g_{jk})\,,
\label{Weyl_zero}\\
\CCC_{ijk} =&\,\nabla_k\RRR_{ij} - \nabla_j\RRR_{ik} -
\frac{1}{4} \big(\nabla_k\RRR g_{ij} - \nabla_j\RRR g_{ik} \big)\,,
\label{Cotton_three}
\end{align}
hence, the evolution equations~\eqref{evolutioncurv} become
\begin{align*}
\frac{\partial\,}{\partial
  t}\Gamma_{ij}^k=&\,-g^{ks}\nabla_i\RRR_{js}-g^{ks}\nabla_j\RRR_{is}+g^{ks}\nabla_s\RRR_{ij}\\
\frac{\partial\,}{\partial t}\RRR_{ij}=&\,\Delta\RRR_{ij}-6g^{pq}\RRR_{ip}\RRR_{jq}
+3\RRR\RRR_{ij}+2\vert\Ric\vert^2g_{ij}-\RRR^2 g_{ij}\\
\frac{\partial\,}{\partial t}\RRR=&\,\Delta\RRR+2|\Ric|^2\,.\\
\end{align*}
From these formulas we can compute the evolution equations of the
derivatives of the curvatures assuming, from now on, to be in normal coordinates,
\begin{eqnarray}
\frac{\partial\,}{\partial
  t}\nabla_l\RRR&=&\nabla_l\Delta\RRR+2\nabla_{l}|\Ric|^2\,,\nonumber\\
&&\nonumber\\
\frac{\partial\,}{\partial t}\nabla_s\RRR_{ij}
&=&\nabla_{s}\Delta\RRR_{ij}-6\nabla_{s}\RRR_{ip}\RRR_{jp}-6\RRR_{ip}\nabla_{s}\RRR_{jp}
+3\nabla_{s}\RRR\RRR_{ij}+3\RRR\nabla_{s}\RRR_{ij}\nonumber\\
&&+2\nabla_{s}\vert\Ric\vert^2g_{ij}-\nabla_{s}\RRR^2 g_{ij}\nonumber\\
&&+(\nabla_i\RRR_{sp}+\nabla_s\RRR_{ip}-\nabla_p\RRR_{is})\RRR_{jp}\nonumber\\
&&+(\nabla_j\RRR_{sp}+\nabla_s\RRR_{jp}-\nabla_p\RRR_{js})\RRR_{ip}\nonumber\\
&=&\nabla_{s}\Delta\RRR_{ij}-5\nabla_{s}\RRR_{ip}\RRR_{jp}-5\RRR_{ip}\nabla_{s}\RRR_{jp}
+3\nabla_{s}\RRR\RRR_{ij}+3\RRR\nabla_{s}\RRR_{ij}\nonumber\\
&&+2\nabla_{s}\vert\Ric\vert^2g_{ij}-\nabla_{s}\RRR^2 g_{ij}\nonumber\\
&&+(\nabla_i\RRR_{sp}-\nabla_p\RRR_{is})\RRR_{jp}+(\nabla_j\RRR_{sp}-\nabla_p\RRR_{js})\RRR_{ip}\nonumber\\
&=&\nabla_{s}\Delta\RRR_{ij}-5\nabla_{s}\RRR_{ip}\RRR_{jp}-5\RRR_{ip}\nabla_{s}\RRR_{jp}
+3\nabla_{s}\RRR\RRR_{ij}+3\RRR\nabla_{s}\RRR_{ij}\nonumber\\
&&+2\nabla_{s}\vert\Ric\vert^2g_{ij}-\nabla_{s}\RRR^2
g_{ij}+\CCC_{spi}\RRR_{jp}+\CCC_{spj}\RRR_{ip}\nonumber\\
&&+\RRR_{jp}[\nabla_i\RRR g_{sp}-\nabla_p\RRR g_{is}]/4
+\RRR_{ip}[\nabla_j\RRR g_{sp}-\nabla_p\RRR g_{js}]/4\,,\nonumber
\end{eqnarray}
where in the last passage we substituted the expression of the Cotton tensor.

We then compute,
\begin{eqnarray*}
\frac{\partial\,}{\partial t}\CCC_{ijk}
&=&\frac{\partial\,}{\partial
  t}\nabla_k\RRR_{ij} -\frac{\partial\,}{\partial t} \nabla_j\RRR_{ik}
-\frac{\partial\,}{\partial t}\big(\nabla_k\RRR g_{ij} - \nabla_j\RRR
g_{ik} \big)/4\\
&=&\nabla_{k}\Delta\RRR_{ij}-5\nabla_{k}\RRR_{ip}\RRR_{jp}-5\RRR_{ip}\nabla_{k}\RRR_{jp}
+3\nabla_{k}\RRR\RRR_{ij}+3\RRR\nabla_{k}\RRR_{ij}\\
&&+2\nabla_{k}\vert\Ric\vert^2g_{ij}-\nabla_{k}\RRR^2
g_{ij}+\CCC_{kpi}\RRR_{jp}+\CCC_{kpj}\RRR_{ip}\\
&&+\RRR_{jp}[\nabla_i\RRR g_{kp}-\nabla_p\RRR g_{ik}]/4
+\RRR_{ip}[\nabla_j\RRR g_{kp}-\nabla_p\RRR g_{jk}]/4\\
&&-\nabla_{j}\Delta\RRR_{ik}+5\nabla_{j}\RRR_{ip}\RRR_{kp}+5\RRR_{ip}\nabla_{j}\RRR_{kp}
-3\nabla_{j}\RRR\RRR_{ik}-3\RRR\nabla_{j}\RRR_{ik}\\
&&-2\nabla_{j}\vert\Ric\vert^2g_{ik}+\nabla_{j}\RRR^2
g_{ik}-\CCC_{jpi}\RRR_{kp}-\CCC_{jpk}\RRR_{ip}\\
&&-\RRR_{kp}[\nabla_i\RRR g_{jp}-\nabla_p\RRR g_{ij}]/4
-\RRR_{ip}[\nabla_k\RRR g_{jp}-\nabla_p\RRR g_{kj}]/4\\
&&+(\RRR_{ij}\nabla_k\RRR -\RRR_{ik}\nabla_j\RRR\big)/2\\
&&-\big(\nabla_k\Delta\RRR+2\nabla_{k}|\Ric|^2\big)g_{ij}/4
+\big(\nabla_j\Delta\RRR+2\nabla_{j}|\Ric|^2\big) g_{ik}/4\\
&=&\nabla_{k}\Delta\RRR_{ij}-5\nabla_{k}\RRR_{ip}\RRR_{jp}-5\RRR_{ip}\nabla_{k}\RRR_{jp}
+3\nabla_{k}\RRR\RRR_{ij}+3\RRR\nabla_{k}\RRR_{ij}\\
&&+3\nabla_{k}\vert\Ric\vert^2g_{ij}/2-\nabla_{k}\RRR^2
g_{ij}+\CCC_{kpi}\RRR_{jp}+\CCC_{kpj}\RRR_{ip}\\
&&+\RRR_{jk}\nabla_i\RRR/4-\RRR_{jp}\nabla_p\RRR g_{ik}/4
+\RRR_{ik}\nabla_j\RRR/4 -\RRR_{ip}\nabla_p\RRR g_{jk}/4\\
&&-\nabla_{j}\Delta\RRR_{ik}+5\nabla_{j}\RRR_{ip}\RRR_{kp}+5\RRR_{ip}\nabla_{j}\RRR_{kp}
-3\nabla_{j}\RRR\RRR_{ik}-3\RRR\nabla_{j}\RRR_{ik}\\
&&-3\nabla_{j}\vert\Ric\vert^2g_{ik}/2+\nabla_{j}\RRR^2
g_{ik}-\CCC_{jpi}\RRR_{kp}-\CCC_{jpk}\RRR_{ip}\\
&&-\RRR_{kj}\nabla_i\RRR/4 +\RRR_{kp}\nabla_p\RRR g_{ij}/4
-\RRR_{ij}\nabla_k\RRR/4+\RRR_{ip}\nabla_p\RRR g_{kj}/4\\
&&+(\RRR_{ij}\nabla_k\RRR -\RRR_{ik}\nabla_j\RRR\big)/2\\
&&-\nabla_k\Delta\RRR g_{ij}/4 + \nabla_j\Delta\RRR g_{ik}/4\\
&=&\nabla_{k}\Delta\RRR_{ij}-5\nabla_{k}\RRR_{ip}\RRR_{jp}-5\RRR_{ip}\nabla_{k}\RRR_{jp}
+13\nabla_{k}\RRR\RRR_{ij}/4+3\RRR\nabla_{k}\RRR_{ij}\\
&&+3\nabla_{k}\vert\Ric\vert^2g_{ij}/2-\nabla_{k}\RRR^2
g_{ij}+\CCC_{kpi}\RRR_{jp}+\CCC_{kpj}\RRR_{ip}\\
&&-\RRR_{jp}\nabla_p\RRR g_{ik}/4\\
&&-\nabla_{j}\Delta\RRR_{ik}+5\nabla_{j}\RRR_{ip}\RRR_{kp}+5\RRR_{ip}\nabla_{j}\RRR_{kp}
-13\nabla_{j}\RRR\RRR_{ik}/4-3\RRR\nabla_{j}\RRR_{ik}\\
&&-3\nabla_{j}\vert\Ric\vert^2g_{ik}/2+\nabla_{j}\RRR^2
g_{ik}-\CCC_{jpi}\RRR_{kp}-\CCC_{jpk}\RRR_{ip}\\
&&+\RRR_{kp}\nabla_p\RRR g_{ij}/4\\
&&-\nabla_k\Delta\RRR g_{ij}/4 + \nabla_j\Delta\RRR g_{ik}/4
\end{eqnarray*}
and
$$
\Delta\CCC_{ijk}=\Delta\nabla_k\RRR_{ij}-\Delta\nabla_j\RRR_{ik}-\Delta\nabla_k\RRR
g_{ij}/4+\Delta\nabla_j\RRR g_{ik}/4\,,
$$
hence,
\begin{eqnarray*}
\frac{\partial\,}{\partial t}\CCC_{ijk}-\Delta\CCC_{ijk}
&=&\nabla_{k}\Delta\RRR_{ij}-\nabla_{j}\Delta\RRR_{ik}-\Delta\nabla_k\RRR_{ij}+\Delta\nabla_j\RRR_{ik}\\
&&-\nabla_k\Delta\RRR g_{ij}/4 + \nabla_j\Delta\RRR
g_{ik}/4+\Delta\nabla_k\RRR g_{ij}/4-\Delta\nabla_j\RRR g_{ik}/4\\
&&-5\nabla_{k}\RRR_{ip}\RRR_{jp}-5\RRR_{ip}\nabla_{k}\RRR_{jp}
+13\nabla_{k}\RRR\RRR_{ij}/4+3\RRR\nabla_{k}\RRR_{ij}\\
&&+3\nabla_{k}\vert\Ric\vert^2g_{ij}/2-\nabla_{k}\RRR^2
g_{ij}+\CCC_{kpi}\RRR_{jp}+\CCC_{kpj}\RRR_{ip}\\
&&-\RRR_{jp}\nabla_p\RRR g_{ik}/4\\
&&+5\nabla_{j}\RRR_{ip}\RRR_{kp}+5\RRR_{ip}\nabla_{j}\RRR_{kp}
-13\nabla_{j}\RRR\RRR_{ik}/4-3\RRR\nabla_{j}\RRR_{ik}\\
&&-3\nabla_{j}\vert\Ric\vert^2g_{ik}/2+\nabla_{j}\RRR^2
g_{ik}-\CCC_{jpi}\RRR_{kp}-\CCC_{jpk}\RRR_{ip}\\
&&+\RRR_{kp}\nabla_p\RRR g_{ij}/4
\end{eqnarray*}

Now to proceed, we need the following commutation rules for the
derivatives of the Ricci tensor and of the scalar curvature, where we will employ the
special form of the Riemann tensor in dimension three given by
formula~\eqref{Weyl_zero},
\begin{eqnarray*}
\nabla_k\Delta\RRR_{ij}-\Delta\nabla_k\RRR_{ij}
&=&\nabla^3_{kll}\RRR_{ij}-\nabla^3_{lkl}\RRR_{ij}+
\nabla^3_{lkl}\RRR_{ij}-\nabla^3_{llk}\RRR_{ij}\\
&=&-\RRR_{kp}\nabla_p\RRR_{ij}
+\RRR_{klip}\nabla_l\RRR_{jp}
+\RRR_{kljp}\nabla_l\RRR_{ip}\\
&&+\nabla^3_{lkl}\RRR_{ij}-\nabla^3_{llk}\RRR_{ij}\\
&=&-\RRR_{kp}\nabla_p\RRR_{ij}
+\RRR_{ik}\nabla_j\RRR/2
+\RRR_{jk}\nabla_i\RRR/2\\
&&-\RRR_{kp}\nabla_i\RRR_{jp}
-\RRR_{kp}\nabla_j\RRR_{ip}
+\RRR_{lp}\nabla_l\RRR_{jp}g_{ik}
+\RRR_{lp}\nabla_l\RRR_{ip}g_{jk}\\
&&-\RRR_{li}\nabla_l\RRR_{jk}
-\RRR_{lj}\nabla_l\RRR_{ik}
-\RRR\nabla_j\RRR g_{ik}/4
-\RRR\nabla_i\RRR g_{jk}/4\\
&&+\RRR\nabla_i\RRR_{jk}/2
+\RRR\nabla_j\RRR_{ik}/2\\
&&+\nabla_l\big(\RRR_{klip}\RRR_{pj}+\RRR_{kljp}\RRR_{pi}\big)\\
&=&-\RRR_{kp}\nabla_p\RRR_{ij}
+\RRR_{ik}\nabla_j\RRR/2
+\RRR_{jk}\nabla_i\RRR/2\\
&&-\RRR_{kp}\nabla_i\RRR_{jp}
-\RRR_{kp}\nabla_j\RRR_{ip}
+\RRR_{lp}\nabla_l\RRR_{jp}g_{ik}
+\RRR_{lp}\nabla_l\RRR_{ip}g_{jk}\\
&&-\RRR_{li}\nabla_l\RRR_{jk}
-\RRR_{lj}\nabla_l\RRR_{ik}
-\RRR\nabla_j\RRR g_{ik}/4
-\RRR\nabla_i\RRR g_{jk}/4\\
&&+\RRR\nabla_i\RRR_{jk}/2
+\RRR\nabla_j\RRR_{ik}/2\\
&&+\nabla_l\big(
\RRR_{ik}\RRR_{lj}
-\RRR_{il}\RRR_{kj}
+\RRR_{pl}\RRR_{pj}g_{ik}
-\RRR_{pk}\RRR_{pj}g_{il}
-g_{ik}\RRR\RRR_{lj}/2
+g_{il}\RRR\RRR_{jk}/2\\
&&
+\RRR_{jk}\RRR_{li}
-\RRR_{jl}\RRR_{ki}
+\RRR_{pl}\RRR_{pi}g_{jk}
-\RRR_{pk}\RRR_{pi}g_{jl}
-g_{jk}\RRR\RRR_{li}/2
+g_{jl}\RRR\RRR_{ik}/2\big)\\
&=&-\RRR_{kp}\nabla_p\RRR_{ij}
+\RRR_{ik}\nabla_j\RRR/2
+\RRR_{jk}\nabla_i\RRR/2\\
&&-\RRR_{kp}\nabla_i\RRR_{jp}
-\RRR_{kp}\nabla_j\RRR_{ip}
+\RRR_{lp}\nabla_l\RRR_{jp}g_{ik}
+\RRR_{lp}\nabla_l\RRR_{ip}g_{jk}\\
&&-\RRR_{li}\nabla_l\RRR_{jk}
-\RRR_{lj}\nabla_l\RRR_{ik}
-\RRR\nabla_j\RRR g_{ik}/4
-\RRR\nabla_i\RRR g_{jk}/4\\
&&+\RRR\nabla_i\RRR_{jk}/2
+\RRR\nabla_j\RRR_{ik}/2\\
&&-\nabla_i\RRR_{pk}\RRR_{pj}
+\nabla_i\RRR\RRR_{jk}/2
+g_{ik}\RRR_{pl}\nabla_l\RRR_{pj}\\
&&-\RRR_{pk}\nabla_i\RRR_{pj}
-g_{ik}\RRR\nabla_j\RRR/4
+\RRR\nabla_i\RRR_{jk}/2\\
&&-\nabla_j\RRR_{pk}\RRR_{pi}
+\nabla_j\RRR\RRR_{ik}/2
+g_{jk}\RRR_{pl}\nabla_l\RRR_{pi}\\
&&-\RRR_{pk}\nabla_j\RRR_{pi}
-g_{jk}\RRR\nabla_i\RRR/4
+\RRR\nabla_j\RRR_{ik}/2\\
&=&-\RRR_{kp}\nabla_p\RRR_{ij}
+\RRR_{ik}\nabla_j\RRR
+\RRR_{jk}\nabla_i\RRR\\
&&-2\RRR_{kp}\nabla_i\RRR_{jp}
-2\RRR_{kp}\nabla_j\RRR_{ip}
+2\RRR_{lp}\nabla_l\RRR_{jp}g_{ik}
+2\RRR_{lp}\nabla_l\RRR_{ip}g_{jk}\\
&&-\RRR_{li}\nabla_l\RRR_{jk}
-\RRR_{lj}\nabla_l\RRR_{ik}
-\RRR_{pj}\nabla_i\RRR_{pk}
-\RRR_{pi}\nabla_j\RRR_{pk}\\
&&-\RRR\nabla_j\RRR g_{ik}/2
-\RRR\nabla_i\RRR g_{jk}/2
+\RRR\nabla_i\RRR_{jk}
+\RRR\nabla_j\RRR_{ik}
\end{eqnarray*}
and 
$$
\nabla_k\Delta\RRR-\Delta\nabla_k\RRR=\RRR_{kllp}\nabla_p\RRR=-\RRR_{kp}\nabla_p\RRR\,.
$$

Then, getting back to the main computation, we obtain
\begin{eqnarray*}
\frac{\partial\,}{\partial t}\CCC_{ijk}-\Delta\CCC_{ijk}
&=&-\RRR_{kp}\nabla_p\RRR_{ij}
+\RRR_{ik}\nabla_j\RRR
+\RRR_{jk}\nabla_i\RRR\\
&&-2\RRR_{kp}\nabla_i\RRR_{jp}
-2\RRR_{kp}\nabla_j\RRR_{ip}
+2\RRR_{lp}\nabla_l\RRR_{jp}g_{ik}
+2\RRR_{lp}\nabla_l\RRR_{ip}g_{jk}\\
&&-\RRR_{li}\nabla_l\RRR_{jk}
-\RRR_{lj}\nabla_l\RRR_{ik}
-\RRR_{pj}\nabla_i\RRR_{pk}
-\RRR_{pi}\nabla_j\RRR_{pk}\\
&&-\RRR\nabla_j\RRR g_{ik}/2
-\RRR\nabla_i\RRR g_{jk}/2
+\RRR\nabla_i\RRR_{jk}
+\RRR\nabla_j\RRR_{ik}\\
&&+\RRR_{jp}\nabla_p\RRR_{ik}
-\RRR_{ij}\nabla_k\RRR
-\RRR_{kj}\nabla_i\RRR\\
&&+2\RRR_{jp}\nabla_i\RRR_{kp}
+2\RRR_{jp}\nabla_k\RRR_{ip}
-2\RRR_{lp}\nabla_l\RRR_{kp}g_{ij}
-2\RRR_{lp}\nabla_l\RRR_{ip}g_{kj}\\
&&+\RRR_{li}\nabla_l\RRR_{kj}
+\RRR_{lk}\nabla_l\RRR_{ij}
+\RRR_{pk}\nabla_i\RRR_{pj}
+\RRR_{pi}\nabla_k\RRR_{pj}\\
&&+\RRR\nabla_k\RRR g_{ij}/2
+\RRR\nabla_i\RRR g_{kj}/2
-\RRR\nabla_i\RRR_{kj}
-\RRR\nabla_k\RRR_{ij}\\
&&+\RRR_{kp}\nabla_p\RRR g_{ij}/4 -\RRR_{jp}\nabla_p\RRR g_{ik}/4\\
&&-5\nabla_{k}\RRR_{ip}\RRR_{jp}-5\RRR_{ip}\nabla_{k}\RRR_{jp}
+13\nabla_{k}\RRR\RRR_{ij}/4+3\RRR\nabla_{k}\RRR_{ij}\\
&&+3\nabla_{k}\vert\Ric\vert^2g_{ij}/2-\nabla_{k}\RRR^2
g_{ij}+\CCC_{kpi}\RRR_{jp}+\CCC_{kpj}\RRR_{ip}\\
&&-\RRR_{jp}\nabla_p\RRR g_{ik}/4\\
&&+5\nabla_{j}\RRR_{ip}\RRR_{kp}+5\RRR_{ip}\nabla_{j}\RRR_{kp}
-13\nabla_{j}\RRR\RRR_{ik}/4-3\RRR\nabla_{j}\RRR_{ik}\\
&&-3\nabla_{j}\vert\Ric\vert^2g_{ik}/2+\nabla_{j}\RRR^2
g_{ik}-\CCC_{jpi}\RRR_{kp}-\CCC_{jpk}\RRR_{ip}\\
&&+\RRR_{kp}\nabla_p\RRR g_{ij}/4\\
&=&\CCC_{kpi}\RRR_{jp}+\CCC_{kpj}\RRR_{ip}-\CCC_{jpi}\RRR_{kp}-\CCC_{jpk}\RRR_{ip}\\
&&+[2\RRR_{lp}\nabla_l\RRR_{jp}+3\RRR\nabla_j\RRR/2
 -\RRR_{jp}\nabla_p\RRR/2-3\nabla_{j}\vert\Ric\vert^2/2]g_{ik}\\
&&+[-2\RRR_{lp}\nabla_l\RRR_{kp}-3\RRR\nabla_k\RRR/2
+\RRR_{kp}\nabla_p\RRR/2+3\nabla_{k}\vert\Ric\vert^2/2]g_{ij}\\
&&-\RRR_{kp}\nabla_i\RRR_{jp}+\RRR_{jp}\nabla_i\RRR_{kp}\\
&&-3\nabla_{k}\RRR_{ip}\RRR_{jp}-4\RRR_{ip}\nabla_k\RRR_{jp}
+9\nabla_{k}\RRR\RRR_{ij}/4+2\RRR\nabla_{k}\RRR_{ij}\\
&&+3\nabla_{j}\RRR_{ip}\RRR_{kp}+4\RRR_{ip}\nabla_{j}\RRR_{kp}
-9\nabla_{j}\RRR\RRR_{ik}/4-2\RRR\nabla_{j}\RRR_{ik}
\end{eqnarray*}

Now, by means of the very definition of the Cotton tensor in
dimension three~\eqref{Cotton_three} and the 
identities~\eqref{CottonSym}, we substitute
\begin{align*}
\CCC_{kpj}-\CCC_{jpk}=&\,-\CCC_{kjp}-\CCC_{jpk}=\CCC_{pkj}\\
\nabla_l\RRR_{jp}=&\,\nabla_j\RRR_{lp}+\CCC_{pjl}+\frac{1}{4}
\big(\nabla_l\RRR g_{pj} - \nabla_j\RRR g_{pl} \big)\\
\nabla_l\RRR_{kp}=&\,\nabla_k\RRR_{lp}+\CCC_{pkl}+\frac{1}{4}
\big(\nabla_l\RRR g_{pk} - \nabla_k\RRR g_{pl} \big)\\
\nabla_i\RRR_{jp}=&\,\nabla_j\RRR_{ip}+\CCC_{pji}+\frac{1}{4}
\big(\nabla_i\RRR g_{jp} - \nabla_j\RRR g_{ip} \big)\\
\nabla_i\RRR_{kp}=&\,\nabla_k\RRR_{ip}+\CCC_{pki}+\frac{1}{4}
\big(\nabla_i\RRR g_{kp} - \nabla_k\RRR g_{ip} \big)
\end{align*}
in the last expression above, getting
\begin{eqnarray*}
\frac{\partial\,}{\partial t}\CCC_{ijk}-\Delta\CCC_{ijk}
&=&\RRR_{jp}\CCC_{kpi}-\RRR_{kp}\CCC_{jpi}+\RRR_{ip}\CCC_{pkj}\\
&&+\Big[2\RRR_{lp}\big(\nabla_j\RRR_{lp}+\CCC_{pjl}+\nabla_l\RRR
g_{pj}/4 - \nabla_j\RRR g_{pl}/4\big)\\
&&\phantom{+\Big[}+3\RRR\nabla_j\RRR/2
 -\RRR_{jp}\nabla_p\RRR/2-3\nabla_{j}\vert\Ric\vert^2/2\Big]g_{ik}\\
&&+\Big[-2\RRR_{lp}\big(\nabla_k\RRR_{lp}+\CCC_{pkl}+\nabla_l\RRR
g_{pk}/4 - \nabla_k\RRR g_{pl}/4\big)\\
&&\phantom{+\Big[}-3\RRR\nabla_k\RRR/2
+\RRR_{kp}\nabla_p\RRR/2+3\nabla_{k}\vert\Ric\vert^2/2\Big]g_{ij}\\
&&-\RRR_{kp}\big(\nabla_j\RRR_{ip}+\CCC_{pji}+\nabla_i\RRR g_{jp}/4 - \nabla_j\RRR g_{ip}/4\big)\\
&&+\RRR_{jp}\big(\nabla_k\RRR_{ip}+\CCC_{pki}+\nabla_i\RRR g_{kp}/4 - \nabla_k\RRR g_{ip}/4\big)\\
&&-3\nabla_{k}\RRR_{ip}\RRR_{jp}-4\RRR_{ip}\nabla_k\RRR_{jp}
+9\nabla_{k}\RRR\RRR_{ij}/4+2\RRR\nabla_{k}\RRR_{ij}\\
&&+3\nabla_{j}\RRR_{ip}\RRR_{kp}+4\RRR_{ip}\nabla_{j}\RRR_{kp}
-9\nabla_{j}\RRR\RRR_{ik}/4-2\RRR\nabla_{j}\RRR_{ik}\\
&=&\RRR_{jp}\big(\CCC_{kpi}+\CCC_{pki}\big)
-\RRR_{kp}\big(\CCC_{jpi}+\CCC_{pji}\big)+\RRR_{ip}\CCC_{pkj}\\
&&+2\RRR_{lp}\CCC_{pjl} g_{ik}-2\RRR_{lp}\CCC_{pkl} g_{ij}\\
&&+\big[\RRR\nabla_j\RRR-\nabla_{j}\vert\Ric\vert^2/2\big]g_{ik}
-\big[\RRR\nabla_k\RRR-\nabla_{k}\vert\Ric\vert^2/2\big]g_{ij}\\
&&-2\nabla_{k}\RRR_{ip}\RRR_{jp}-4\RRR_{ip}\nabla_k\RRR_{jp}
+2\nabla_{k}\RRR\RRR_{ij}+2\RRR\nabla_{k}\RRR_{ij}\\
&&+2\nabla_{j}\RRR_{ip}\RRR_{kp}+4\RRR_{ip}\nabla_{j}\RRR_{kp}
-2\nabla_{j}\RRR\RRR_{ik}-2\RRR\nabla_{j}\RRR_{ik}\,.
\end{eqnarray*}
then, we substitute again
\begin{align*}
\nabla_k\RRR_{jp}=&\,\nabla_p\RRR_{kj}+\CCC_{jpk}+\frac{1}{4}
\big(\nabla_k\RRR g_{jp} - \nabla_p\RRR g_{jk} \big)\\
\nabla_j\RRR_{kp}=&\,\nabla_p\RRR_{jk}+\CCC_{kpj}+\frac{1}{4}
\big(\nabla_j\RRR g_{kp} - \nabla_p\RRR g_{kj} \big)\\
\nabla_k\RRR_{ij}=&\,\nabla_i\RRR_{kj}+\CCC_{jik}+\frac{1}{4}
\big(\nabla_k\RRR g_{ij} - \nabla_i\RRR g_{jk} \big)\\
\nabla_j\RRR_{ik}=&\,\nabla_i\RRR_{jk}+\CCC_{kij}+\frac{1}{4}
\big(\nabla_j\RRR g_{ik} - \nabla_i\RRR g_{kj} \big)\,,
\end{align*}
finally obtaining
\begin{eqnarray*}
\frac{\partial\,}{\partial t}\CCC_{ijk}-\Delta\CCC_{ijk}
&=&\RRR_{jp}\big(\CCC_{kpi}+\CCC_{pki}\big)
-\RRR_{kp}\big(\CCC_{jpi}+\CCC_{pji}\big)+\RRR_{ip}\CCC_{pkj}\\
&&+2\RRR_{lp}\CCC_{pjl} g_{ik}-2\RRR_{lp}\CCC_{pkl} g_{ij}\\
&&+\big[\RRR\nabla_j\RRR-\nabla_{j}\vert\Ric\vert^2/2\big]g_{ik}
-\big[\RRR\nabla_k\RRR-\nabla_{k}\vert\Ric\vert^2/2\big]g_{ij}\\
&&-2\nabla_{k}\RRR_{ip}\RRR_{jp}-4\RRR_{ip}\big(
\nabla_p\RRR_{kj}+\CCC_{jpk}+\nabla_k\RRR g_{jp}/4 - \nabla_p\RRR g_{jk}/4 \big)\\
&&+2\nabla_{k}\RRR\RRR_{ij}
+2\RRR\big(
\nabla_i\RRR_{kj}+\CCC_{jik}+\nabla_k\RRR g_{ij}/4 - \nabla_i\RRR g_{jk}/4 \big)\\
&&+2\nabla_{j}\RRR_{ip}\RRR_{kp}+4\RRR_{ip}\big(
\nabla_p\RRR_{jk}+\CCC_{kpj}+\nabla_j\RRR g_{kp}/4 - \nabla_p\RRR g_{kj}/4 \big)\\
&&-2\nabla_{j}\RRR\RRR_{ik}
-2\RRR\big(
\nabla_i\RRR_{jk}+\CCC_{kij}+\nabla_j\RRR g_{ik}/4 - \nabla_i\RRR
g_{kj}/4 \big)\\
&=&\RRR_{jp}\big(\CCC_{kpi}+\CCC_{pki}\big)
-\RRR_{kp}\big(\CCC_{jpi}+\CCC_{pji}\big)+\RRR_{ip}\CCC_{pkj}\\
&&+4\RRR_{ip}\big(\CCC_{kpj}-\CCC_{jpk}\big)
+2\RRR\big(\CCC_{jik}-\CCC_{kij}\big)\\
&&+2\RRR_{lp}\CCC_{pjl} g_{ik}-2\RRR_{lp}\CCC_{pkl} g_{ij}\\
&&+\big[\RRR\nabla_j\RRR/2-\nabla_{j}\vert\Ric\vert^2/2\big]g_{ik}
-\big[\RRR\nabla_k\RRR/2-\nabla_{k}\vert\Ric\vert^2/2\big]g_{ij}\\
&&-2\nabla_{k}\RRR_{ip}\RRR_{jp}+2\nabla_{j}\RRR_{ip}\RRR_{kp}\\
&&+\nabla_{k}\RRR\RRR_{ij}-\nabla_{j}\RRR\RRR_{ik}\\
&=&\RRR_{jp}\big(\CCC_{kpi}+\CCC_{pki}\big)
-\RRR_{kp}\big(\CCC_{jpi}+\CCC_{pji}\big)+5\RRR_{ip}\CCC_{pkj}\\
&&+2\RRR\CCC_{ijk}+2\RRR_{lp}\CCC_{pjl} g_{ik}-2\RRR_{lp}\CCC_{pkl} g_{ij}\\
&&+\big[\RRR\nabla_j\RRR/2-\nabla_{j}\vert\Ric\vert^2/2\big]g_{ik}
-\big[\RRR\nabla_k\RRR/2-\nabla_{k}\vert\Ric\vert^2/2\big]g_{ij}\\
&&+2\nabla_{j}\RRR_{ip}\RRR_{kp}-2\nabla_{k}\RRR_{ip}\RRR_{jp}\\
&&+\nabla_{k}\RRR\RRR_{ij}-\nabla_{j}\RRR\RRR_{ik}\,,
\end{eqnarray*}
where in the last passage we used again the 
identities~\eqref{CottonSym}.\\
Hence, we can resume this long computation in the following
proposition, getting back to a generic coordinate basis.
\begin{prop}\label{cot}
During the Ricci flow of a 3--dimensional Riemannian manifold $(M^3,
g(t))$, 
the Cotton tensor satisfies the following evolution equation
\begin{eqnarray}\label{main_eq1}
\bigl(\partial_t-\Delta\bigr)\CCC_{ijk}&=&g^{pq}\RRR_{pj}(\CCC_{kqi}+\CCC_{qki})+5g^{pq}\RRR_{ip}\CCC_{qkj}
+g^{pq}\RRR_{pk}(\CCC_{jiq}+\CCC_{qij})\\
&&+2\RRR\CCC_{ijk}+2\RRR^{ql}\CCC_{qjl}g_{ik}-2\RRR^{ql}\CCC_{qkl}g_{ij}\nonumber\\
&&+\frac{1}{2}\nabla_k|\Ric|^2g_{ij}-\frac{1}{2}\nabla_j|\Ric|^2g_{ik}+\frac{\RRR}{2}\nabla_j\RRR g_{ik}
-\frac{\RRR}{2}\nabla_k\RRR g_{ij}\nonumber\\
&&+2g^{pq}\RRR_{pk}\nabla_j\RRR_{qi}-2g^{pq}\RRR_{pj}\nabla_k\RRR_{qi}
+\RRR_{ij}\nabla_k\RRR-\RRR_{ik}\nabla_j\RRR\,.\nonumber
\end{eqnarray}
In particular if the Cotton tensor vanishes identically along the flow we obtain,
\begin{eqnarray}\label{main_eq2}
0&=&\nabla_k|\Ric|^2g_{ij}-\nabla_j|\Ric|^2g_{ik}+{\RRR}\nabla_j\RRR g_{ik}
-{\RRR}\nabla_k\RRR g_{ij}\\
&&+4g^{pq}\RRR_{pk}\nabla_j\RRR_{qi}-4g^{pq}\RRR_{pj}\nabla_k\RRR_{qi}
+2\RRR_{ij}\nabla_k\RRR-2\RRR_{ik}\nabla_j\RRR\,.\nonumber
\end{eqnarray}
\end{prop}

\begin{cor}
If the Cotton tensor vanishes identically along the Ricci flow of a 3--dimensional Riemannian manifold $(M^3, g(t))$, the following tensor
$$
|\Ric|^2g_{ij}-4\RRR_{pj}\RRR_{pi}+3\RRR\RRR_{ij}-\frac{7}{8}\RRR^2g_{ij}
$$
is a Codazzi tensor (see~\cite[Chapter~16, Section~C]{besse}).
\end{cor}
\begin{proof}
We compute in an orthonormal basis,
\begin{align*}
4\RRR_{pk}\nabla_j\RRR_{pi}-4&\,\RRR_{pj}\nabla_k\RRR_{pi}+2\RRR_{ij}\nabla_k\RRR-2\RRR_{ik}\nabla_j\RRR\\
=&\,4\nabla_j(\RRR_{pk}\RRR_{pi})-4\nabla_k(\RRR_{pj}\RRR_{pi})
-4\RRR_{pi}\nabla_j\RRR_{pk}+4\RRR_{pi}\nabla_k\RRR_{pj}\\
&\,+2\RRR_{ij}\nabla_k\RRR-2\RRR_{ik}\nabla_j\RRR\\
=&\,4\nabla_j(\RRR_{pk}\RRR_{pi})-4\nabla_k(\RRR_{pj}\RRR_{pi})
+\RRR_{pi}(4\CCC_{pjk}+\nabla_k\RRR g_{pj}-\nabla_j\RRR g_{pk})\\
&\,+2\RRR_{ij}\nabla_k\RRR-2\RRR_{ik}\nabla_j\RRR\\
=&\,4\nabla_j(\RRR_{pk}\RRR_{pi})-4\nabla_k(\RRR_{pj}\RRR_{pi})
+3\RRR_{ij}\nabla_k\RRR-3\RRR_{ik}\nabla_j\RRR\\
=&\,4\nabla_j(\RRR_{pk}\RRR_{pi})-4\nabla_k(\RRR_{pj}\RRR_{pi})
+3\nabla_k(\RRR\RRR_{ij})-3\nabla_j(\RRR\RRR_{ik})\\
&\,-3\RRR(\nabla_k\RRR_{ij}-\nabla_j\RRR_{ik})\\
=&\,4\nabla_j(\RRR_{pk}\RRR_{pi})-4\nabla_k(\RRR_{pj}\RRR_{pi})
+3\nabla_k(\RRR\RRR_{ij})-3\nabla_j(\RRR\RRR_{ik})\\
&\,-3\RRR(4\CCC_{ijk}+\nabla_k\RRR g_{ij}-\nabla_j\RRR g_{ik})/4\\
=&\,4\nabla_j(\RRR_{pk}\RRR_{pi})-4\nabla_k(\RRR_{pj}\RRR_{pi})
+3\nabla_k(\RRR\RRR_{ij})-3\nabla_j(\RRR\RRR_{ik})\\
&\,-\frac{3}{8}\nabla_k\RRR^2 g_{ij}+\frac{3}{8}\nabla_j\RRR^2 g_{ik}\,.
\end{align*}
Hence, we have, by the previous proposition,
\begin{align*}
0=&\,\nabla_k|\Ric|^2g_{ij}-\nabla_j|\Ric|^2g_{ik}
+4\nabla_j(\RRR_{pk}\RRR_{pi})-4\nabla_k(\RRR_{pj}\RRR_{pi})\\
&\,+3\nabla_k(\RRR\RRR_{ij})-3\nabla_j(\RRR\RRR_{ik})
-\frac{7}{8}\nabla_k\RRR^2 g_{ij}+\frac{7}{8}\nabla_j\RRR^2 g_{ik}\,,
\end{align*}
which is the thesis of the corollary.
\end{proof}

\begin{rem}
All the traces of the 3--tensor in the LHS of equation~\eqref{main_eq2} are zero.
\end{rem}

\begin{rem}
From the trace--free property~\eqref{CottonTraces} of the Cotton
tensor and the fact that along the Ricci flow there holds
$$
\bigl(\partial_t-\Delta\bigr)g^{ij}=2\RRR^{ij}\,,
$$
we conclude that the following relations have to hold
\begin{eqnarray*}
g^{ij}(\partial_t-\Delta) \CCC_{ijk}&=&-2\RRR^{ij}\CCC_{ijk}\,,\\
g^{ik}(\partial_t-\Delta) \CCC_{ijk}&=&-2\RRR^{ik}\CCC_{ijk}\,,\\
g^{jk}(\partial_t-\Delta) \CCC_{ijk}&=&-2\RRR^{jk}\CCC_{ijk}=0\,.\
\end{eqnarray*}
They are easily verified for formula~\eqref{main_eq1}.
\end{rem}

\begin{cor}\label{CorEv}
During the Ricci flow of a 3--dimensional Riemannian manifold $(M^3,
g(t))$, the squared norm of the Cotton tensor satisfies the following evolution equation, in an orthonormal basis,
\begin{eqnarray*}
\bigl(\partial_t-\Delta\bigr)\vert\CCC_{ijk}\vert^2
&=&-2\vert \nabla \CCC_{ijk}\vert^2-16\CCC_{ipk}\CCC_{iqk}\RRR_{pq}
+24\CCC_{ipk}\CCC_{kqi}\RRR_{pq}+4\RRR\vert\CCC_{ijk}\vert^2\\
&&+8\CCC_{ijk}\RRR_{pk}\nabla_j\RRR_{pi}
+4\CCC_{ijk}\RRR_{ij}\nabla_k\RRR\,.\nonumber
\end{eqnarray*}
\end{cor}
\begin{proof}
\begin{eqnarray*}
\bigl(\partial_t-\Delta\bigr)\vert\CCC_{ijk}\vert^2
&=&-2\vert \nabla \CCC_{ijk}\vert^2
+2\CCC^{ijk}\RRR_{ip}g^{pq}\CCC_{qjk}
+2\CCC^{ijk}\RRR_{jp}g^{pq}\CCC_{iqk}
+2\CCC^{ijk}\RRR_{kp}g^{pq}\CCC_{ijq}\\
&&+2\CCC^{ijk}\Bigl[g^{pq}\RRR_{pj}(\CCC_{kqi}+\CCC_{qki})+5g^{pq}\RRR_{ip}\CCC_{qkj}
+g^{pq}\RRR_{pk}(\CCC_{jiq}+\CCC_{qij})\\
&&+2\RRR\CCC_{ijk}+2\RRR^{ql}\CCC_{qjl}g_{ik}-2\RRR^{ql}\CCC_{qkl}g_{ij}\\
&&+\frac{1}{2}\nabla_k|\Ric|^2g_{ij}-\frac{1}{2}\nabla_j|\Ric|^2g_{ik}+\frac{\RRR}{2}\nabla_j\RRR g_{ik}
-\frac{\RRR}{2}\nabla_k\RRR g_{ij}\\
&&+2g^{pq}\RRR_{pk}\nabla_j\RRR_{qi}-2g^{pq}\RRR_{pj}\nabla_k\RRR_{qi}
+\RRR_{ij}\nabla_k\RRR-\RRR_{ik}\nabla_j\RRR\Bigr]\\
&=&-2\vert \nabla \CCC_{ijk}\vert^2
+2(\CCC^{kij}+\CCC^{jki})\RRR_{ip}g^{pq}(\CCC_{kqj}+\CCC_{jkq})\\
&&+2\CCC^{ijk}\RRR_{jp}g^{pq}\CCC_{iqk}
+2\CCC^{ikj}\RRR_{kp}g^{pq}\CCC_{iqj}\\
&&+2\CCC^{ijk}\Bigl[2g^{pq}\RRR_{pj}(\CCC_{kqi}+\CCC_{qki})+5g^{pq}\RRR_{ip}\CCC_{qkj}\Bigr]\\
&&+4\RRR\vert\CCC_{ijk}\vert^2
+8g^{pq}\CCC^{ijk}\RRR_{pk}\nabla_j\RRR_{qi}
+4\CCC^{ijk}\RRR_{ij}\nabla_k\RRR\\
&=&-2\vert \nabla \CCC_{ijk}\vert^2-16\CCC_{ipk}\CCC_{iqk}\RRR_{pq}
+24\CCC_{ipk}\CCC_{kqi}\RRR_{pq}+4\RRR\vert\CCC_{ijk}\vert^2\\
&&+8\CCC_{ijk}\RRR_{pk}\nabla_j\RRR_{pi}
+4\CCC_{ijk}\RRR_{ij}\nabla_k\RRR
\end{eqnarray*}
where in the last line we assumed to be in a orthonormal basis.
\end{proof}

\section{Three--Dimensional Gradient Ricci Solitons}\label{cgrad}

The structural equation of a gradient Ricci soliton $(M^n,g, \nabla f)$ is the following
\begin{equation}\label{SolEq0}
\RRR_{ij}+\nabla_i\nabla_jf=\lambda g_{ij}\,,
\end{equation}
for some $\lambda\in\mathbb{R}$.\\
The soliton is said to be {\em steady}, {\em shrinking} or {\em expanding} according to the fact that the constant $\lambda$ is zero, positive or negative, respectively.

It follows that in dimension three, for $(M^3,g, \nabla f)$ there holds
\begin{eqnarray}
\Delta\RRR_{ij}&=&\nabla_l\RRR_{ij}\nabla_l
f+2\lambda\RRR_{ij}-2|\Ric|^2g_{ij}+\RRR^2
g_{ij}-3\RRR\RRR_{ij}+4\RRR_{is}\RRR_{sj}\,\label{SolEq1}\\
\Delta\RRR&=&\nabla_l\RRR\nabla_l f+2\lambda\RRR-2|\Ric|^2\,\label{SolEq2}\\
\nabla_i\RRR&=&2\RRR_{li}\nabla_l f\label{SolEq3}\\
\CCC_{ijk}&=&\frac{\RRR_{lk}g_{ij}}{2}\nabla_lf-\frac{\RRR_{lj}g_{ik}}{2}\nabla_l
f+\RRR_{ij}\nabla_k f-\RRR_{ik}\nabla_j f
+\frac{\RRR g_{ik}}{2}\nabla_j f -\frac{\RRR g_{ij}}{2}\nabla_k f\label{SolEq4}\\
&=&\frac{\nabla_k\RRR}{4}g_{ij}-\frac{\nabla_j\RRR}{4}g_{ik}
+\Bigl(\RRR_{ij}-\frac{\RRR}{2}g_{ij}\Bigr)\nabla_k f
-\Bigl(\RRR_{ik}-\frac{\RRR}{2}g_{ik}\Bigr)\nabla_j f\,.\nonumber
\end{eqnarray}
In the special case of a {\em steady} soliton the first two equations above simplify as follows,
\begin{eqnarray*}
\Delta\RRR_{ij}&=&\nabla_l\RRR_{ij}\nabla_l
f-2|\Ric|^2g_{ij}+\RRR^2
g_{ij}-3\RRR\RRR_{ij}+4\RRR_{is}\RRR_{sj}\\
\Delta\RRR&=&\nabla_l\RRR\nabla_l f-2|\Ric|^2\,.
\end{eqnarray*}

\begin{rem}
We notice that, by relation~\eqref{SolEq4}, we have
\begin{eqnarray*}
\CCC_{ijk}\nabla_i f&=&\frac{\nabla_k\RRR\nabla_jf}{4}-\frac{\nabla_j\RRR\nabla_k f}{4}+\RRR_{ij}\nabla_if\nabla_kf-\frac{\RRR}{2}\nabla_jf\nabla_k f
-\RRR_{ik}\nabla_if\nabla_jf+\frac{\RRR}{2}\nabla_kf\nabla_j f\\
&=&\frac{\nabla_j\RRR\nabla_kf}{4}-\frac{\nabla_k\RRR\nabla_j f}{4}\,,
\end{eqnarray*}
where in the last passage we used relation~\eqref{SolEq3}.\\
It follows that
$$
\CCC_{ijk}\nabla_i f\nabla_j f=\frac{\langle\nabla f,\nabla\RRR\rangle}{4}\nabla_k f-\frac{\vert\nabla f\vert^2}{4}\nabla_k\RRR\,.
$$
Hence, if the Cotton tensor of a three--dimensional 
gradient Ricci soliton is identically zero, we have that at every point where $\nabla\RRR$ is not zero, $\nabla f$ and $\nabla\RRR$ are proportional.

This relation is a key step in (yet another) proof of the fact that a
three--dimensional, locally conformally flat, steady or
shrinking gradient Ricci soliton is locally a warped product of a 
constant curvature surface on a
interval of $\RR$, leading to a full classification, first obtained by
H.-D. Cao and Q. Chen~\cite{caochen} for the steady case and
H.-D. Cao, B.-L. Chen and X.-P. Zhu~\cite{caochenzhu} for the
shrinking case (actually this is the last paper of a series finally
classifying, in full generality, all the three-dimensional gradient
shrinking Ricci solitons, even without the LCF assumption).
\end{rem}

\begin{prop}
Let $(M^3, g, f)$ be a three--dimensional gradient Ricci soliton. Then, 
\begin{eqnarray*}
\Delta|\CCC_{ijk}|^2&=&\nabla_l|\CCC_{ijk}|^2
\nabla_l f+2|\nabla\CCC_{ijk}|^2-2\RRR|\CCC_{ijk}|^2\\
&&-6\CCC_{ijk}\RRR_{ij}\nabla_k\RRR+8\CCC_{jsk}\CCC_{jik}\RRR_{si}-16\CCC_{jsk}\CCC_{kij}\RRR_{si}
-8\CCC_{ijk}\RRR_{lk}\nabla_j\RRR_{il}\,.
\end{eqnarray*}
\end{prop}
\begin{proof}
First observe that
\[
\ \Delta|\CCC_{ijk}|^2
=2\CCC_{ijk}\Delta\CCC_{ijk}+2|\nabla\CCC_{ijk}|^2.
\]
Using relations~\eqref{SolEq4},~\eqref{SolEq1} and, repeatedly, the trace--free property~\eqref{CottonTraces} of the Cotton tensor, we have that
\begin{eqnarray*}
\CCC_{ijk}\Delta\CCC_{ijk}&=&\Delta(\RRR_{ij}\nabla_k f-\RRR_{ik}\nabla_j f)\CCC_{ijk}\\
&=&(\Delta\RRR_{ij}\nabla_k f+\RRR_{ij}\Delta\nabla_k f+2\nabla_l\RRR_{ij}\nabla_l\nabla_k f)\CCC_{ijk}\\
&&-(\Delta\RRR_{ik}\nabla_j f+\RRR_{ik}\Delta\nabla_j f+2\nabla_l\RRR_{ik}\nabla_l\nabla_j f)\CCC_{ijk}\\
&=&(\nabla_s\RRR_{ij}\nabla_s f-3\RRR\RRR_{ij}+4\RRR_{is}\RRR_{sj})\nabla_k f\CCC_{ijk}\\
&&+\RRR_{ij}\Delta\nabla_k f\CCC_{ijk}+2\nabla_l\RRR_{ij}\nabla_l\nabla_k f \CCC_{ijk}\\
&&-(\nabla_s\RRR_{ik}\nabla_s f-3\RRR\RRR_{ik}+4\RRR_{is}\RRR_{sk})\nabla_j f\CCC_{ijk}\\
&&-\RRR_{ik}\Delta\nabla_j f\CCC_{ijk}-2\nabla_l\RRR_{ik}\nabla_l\nabla_j f \CCC_{ijk}\\
&=&(\nabla_s\RRR_{ij}\nabla_k f- \nabla_s\RRR_{ik}\nabla_j f)\nabla_s f\CCC_{ijk}\\
&&-3\RRR(\RRR_{ij}\nabla_k f-\RRR_{ik}\nabla_j f)\CCC_{ijk}\\
&&+4\RRR_{is}(\RRR_{sj}\nabla_k f-\RRR_{sk}\nabla_j f)\CCC_{ijk}\\
&&+(\RRR_{ij}\nabla_l\nabla_l\nabla_k f-\RRR_{ik}\nabla_l\nabla_l\nabla_j f)\CCC_{ijk}\\
&&+2(\nabla_l\RRR_{ij}\nabla_l\nabla_k f-\nabla_l\RRR_{ik}\nabla_l\nabla_j f)\CCC_{ijk}\\
&=&(\nabla_s\RRR_{ij}\nabla_k f- \nabla_s\RRR_{ik}\nabla_j f)\nabla_s f\CCC_{ijk}\\
&&+(-3\RRR)\vert\CCC_{ijk}\vert^2\\
&&+4\RRR_{is}(\RRR_{sj}\nabla_k f-\RRR_{sk}\nabla_j f)\CCC_{ijk}\\
&&+(\RRR_{ij}\nabla_l\nabla_l\nabla_k f-\RRR_{ik}\nabla_l\nabla_l\nabla_j f)\CCC_{ijk}\\
&&+2(\nabla_l\RRR_{ij}\nabla_l\nabla_k f-\nabla_l\RRR_{ik}\nabla_l\nabla_j f)\CCC_{ijk}\,,
\end{eqnarray*}
where we used the identity
\begin{equation}\label{eq_NormaCottonSol}
(\RRR_{ij}\nabla_k f-\RRR_{ik}\nabla_j f)\CCC_{ijk}=\vert\CCC_{ijk}\vert^2
\end{equation}
which follows easily by equation~\eqref{SolEq4} and the fact that every trace of the Cotton tensor is zero.\\
Using now equations~\eqref{SolEq0},~\eqref{SolEq4},~\eqref{CottonTraces},~\eqref{CottonSym}, and~\eqref{SolEq3}, we compute
\begin{eqnarray*}
(\nabla_s\RRR_{ij}\nabla_k f-\nabla_s\RRR_{ik}\nabla_j f)\nabla_s f\CCC_{ijk}&=&(\nabla_s(\RRR_{ij}\nabla_k f)-\RRR_{ij}\nabla_s\nabla_k f)\nabla_s f\CCC_{ijk}\\
&&-(\nabla_s(\RRR_{ik}\nabla_j f)-\RRR_{ik}\nabla_s\nabla_j f)\nabla_s f\CCC_{ijk}\\
&=&(\nabla_s(\RRR_{ij}\nabla_k f-\RRR_{ik}\nabla_j f)+\RRR_{ij}(\RRR_{sk}))\nabla_s f\CCC_{ijk}\\
&&-(\RRR_{ik}(\RRR_{sj}))\nabla_s f\CCC_{ijk}\\
&=&\nabla_s\CCC_{ijk}\CCC_{ijk}\nabla_s f+\RRR_{ij}\RRR_{sk}\nabla_s f\CCC_{ijk}-\RRR_{ik}\RRR_{sj}\nabla_s f \CCC_{ijk}\\
&=&\frac{1}{2}\nabla_s|\CCC_{ijk}|^2\nabla_s f+\frac{1}{2}\RRR_{ij}\nabla_k \RRR\CCC_{ijk}-\frac{1}{2}\RRR_{ik}\nabla_j \RRR\CCC_{ijk}\\
&&\\
&&\\
4\RRR_{is}(\RRR_{sj}\nabla_k f-\RRR_{sk}\nabla_j f)\CCC_{ijk}
&=&4\RRR_{is}(\CCC_{sjk}-\frac{1}{4}\nabla_k\RRR g_{sj}+\frac{1}{4}\nabla_j\RRR g_{sk}+\frac{\RRR}{2}\nabla_k f g_{sj}-\frac{\RRR}{2}\nabla_j f g_{sk})\CCC_{ijk}\\
&=&4\RRR_{is}(-\CCC_{jks}-\CCC_{ksj})(-\CCC_{jki}-\CCC_{kij})-\RRR_{ij}\nabla_k\RRR\CCC_{ijk}\\
&&+\RRR_{ik}\nabla_j\RRR\CCC_{ijk}+2\RRR\RRR_{ij}\nabla_k f\CCC_{ijk}-2\RRR\RRR_{ik}\nabla_j f\CCC_{ijk}\\
&=&8\RRR_{is}\CCC_{jsk}\CCC_{jik}-8\RRR_{is}\CCC_{jsk}\CCC_{kij}\\
&&-\RRR_{ij}\nabla_k\RRR\CCC_{ijk}+\RRR_{ik}\nabla_j\RRR\CCC_{ijk}+2\RRR|\CCC_{ijk}|^2\\
&&\\
&&\\
(\RRR_{ij}\nabla_l\nabla_l\nabla_k f-\RRR_{ik}\nabla_l\nabla_l\nabla_j
f)\CCC_{ijk}&=&(\RRR_{ij}\nabla_l(-\RRR_{lk})
-\RRR_{ik}\nabla_l(-\RRR_{lj}))\CCC_{ijk}\\
&=&-\frac{1}{2}\RRR_{ij}\nabla_k\RRR\CCC_{ijk}+\frac{1}{2}\RRR_{ik}\nabla_j\RRR\CCC_{ijk}\\
&&\\
&&\\
2(\nabla_l\RRR_{ij}\nabla_l\nabla_k f-\nabla_l\RRR_{ik}\nabla_l\nabla_j f)\CCC_{ijk}&=&2((\CCC_{ijl}+\nabla_j\RRR_{il}+\frac{1}{4}\nabla_l\RRR g_{ij}-\frac{1}{4}\nabla_j\RRR g_{il})(-\RRR_{lk}))\CCC_{ijk}\\
&&-2((\CCC_{ikl}+\nabla_k\RRR_{il}+\frac{1}{4}\nabla_l\RRR g_{ik}-\frac{1}{4}\nabla_k\RRR g_{il})(-\RRR_{lj}))\CCC_{ijk}\\
&=&-2\CCC_{ijl}\CCC_{ijk}\RRR_{lk}-2\CCC_{ijk}\RRR_{lk}\nabla_j\RRR_{il}+\frac{1}{2}\CCC_{ijk}\RRR_{ik}\nabla_j\RRR\\
&&+2\CCC_{ikl}\CCC_{ijk}\RRR_{lj}+2\CCC_{ijk}\RRR_{lj}\nabla_k\RRR_{il}-\frac{1}{2}\CCC_{ijk}\RRR_{ij}\nabla_k\RRR\\
&=&-2\CCC_{ilj}\CCC_{ikj}\RRR_{lk}-2\CCC_{ijk}\RRR_{lk}\nabla_j\RRR_{il}+\frac{1}{2}\CCC_{ijk}\RRR_{ik}\nabla_j\RRR\\
&&-2\CCC_{ilk}\CCC_{ijk}\RRR_{lj}+2\CCC_{ijk}\RRR_{lj}\nabla_k\RRR_{il}-\frac{1}{2}\CCC_{ijk}\RRR_{ij}\nabla_k\RRR.
\end{eqnarray*}
Hence, getting back to the main computation and using again the symmetry relations~\eqref{CottonSym}, we finally get
\begin{eqnarray*}
\CCC_{ijk}\Delta\CCC_{ijk}&=&\frac{1}{2}\nabla_s|\CCC_{ijk}|^2\nabla_s f-\RRR|\CCC_{ijk}|^2\\
&&-\frac{3}{2}\CCC_{ijk}\RRR_{ij}\nabla_k\RRR+\frac{3}{2}\CCC_{ijk}\RRR_{ik}\nabla_j\RRR\\
&&+4\CCC_{jsk}\CCC_{jik}\RRR_{si}-8\CCC_{jsk}\CCC_{kij}\RRR_{si}\\
&&-2\CCC_{ijk}\RRR_{lk}\nabla_j\RRR_{il}+2\CCC_{ijk}\RRR_{lj}\nabla_k\RRR_{il}\\
&=&\frac{1}{2}\nabla_s|\CCC_{ijk}|^2\nabla_s f-\RRR|\CCC_{ijk}|^2\\
&&-3\CCC_{ijk}\RRR_{ij}\nabla_k\RRR+4\CCC_{jsk}\CCC_{jik}\RRR_{si}-8\CCC_{jsk}\CCC_{kij}\RRR_{si}-4\CCC_{ijk}\RRR_{lk}\nabla_j\RRR_{il}\,
\end{eqnarray*}
where in the last passage we applied the skew--symmetry of the Cotton tensor in its last two indexes. The thesis follows.
\end{proof}

\section{The Evolution Equation of the Cotton Tensor in any Dimension}\label{nsec}

In this section we will compute the evolution equation under
the Ricci flow of the Cotton tensor $\CCC_{ijk}$, for every
$n$--dimensional Riemannian manifold $(M^n,g(t))$ evolving by Ricci flow.

Among the evolution equations~\eqref{evolutioncurv} we expand the one for the Ricci tensor,
\begin{align*}
\frac{\partial\,}{\partial t}\RRR_{ij}=&\,\Delta\RRR_{ij}-\frac{2n}{n-2}g^{pq}\RRR_{ip}\RRR_{jq}
+\frac{2n}{(n-1)(n-2)}\RRR\RRR_{ij}+\frac{2}{n-2}\vert\Ric\vert^2g_{ij}\\
&\,-\frac{2}{(n-1)(n-2)}\RRR^2 g_{ij}-2\RRR^{pq}\WWW_{pijq}\,.
\end{align*}
Then, we compute the evolution equations of the
derivatives of the curvatures assuming, from now on, to be in normal coordinates,
\begin{eqnarray}
\frac{\partial\,}{\partial
  t}\nabla_l\RRR&=&\nabla_l\Delta\RRR+2\nabla_{l}|\Ric|^2\,,\nonumber\\
&&\nonumber\\
\frac{\partial\,}{\partial t}\nabla_s\RRR_{ij}
&=&\nabla_{s}\Delta\RRR_{ij}-\frac{2n}{n-2}\nabla_{s}\RRR_{ip}\RRR_{jp}-\frac{2n}{n-2}\RRR_{ip}\nabla_{s}\RRR_{jp}
+\frac{2n}{(n-1)(n-2)}\nabla_{s}\RRR\RRR_{ij}\nonumber\\
&&+\frac{2n}{(n-1)(n-2)}\RRR\nabla_{s}\RRR_{ij}+\frac{2}{n-2}\nabla_{s}\vert\Ric\vert^2g_{ij}-\frac{2}{(n-1)(n-2)}\nabla_{s}\RRR^2 g_{ij}\nonumber\\
&&-2\nabla_{s}\RRR_{kl}\WWW_{kijl}-2\RRR_{kl}\nabla_{s}\WWW_{kijl}+(\nabla_i\RRR_{sp}+\nabla_s\RRR_{ip}-\nabla_p\RRR_{is})\RRR_{jp}\nonumber\\
&&+(\nabla_j\RRR_{sp}+\nabla_s\RRR_{jp}-\nabla_p\RRR_{js})\RRR_{ip}\nonumber\\
&=&\nabla_{s}\Delta\RRR_{ij}-\frac{n+2}{n-2}\nabla_{s}\RRR_{ip}\RRR_{jp}-\frac{n+2}{n-2}\RRR_{ip}\nabla_{s}\RRR_{jp}
+\frac{2n}{(n-1)(n-2)}\nabla_{s}\RRR\RRR_{ij}\nonumber\\
&&+\frac{2n}{(n-1)(n-2)}\RRR\nabla_{s}\RRR_{ij}+\frac{2}{n-2}\nabla_{s}\vert\Ric\vert^2g_{ij}-\frac{2}{(n-1)(n-2)}\nabla_{s}\RRR^2 g_{ij}\nonumber\\
&&-2\nabla_{s}\RRR_{kl}\WWW_{kijl}-2\RRR_{kl}\nabla_{s}\WWW_{kijl}+(\nabla_i\RRR_{sp}-\nabla_p\RRR_{is})\RRR_{jp}+(\nabla_j\RRR_{sp}-\nabla_p\RRR_{js})\RRR_{ip}\nonumber\\
&=&\nabla_{s}\Delta\RRR_{ij}-\frac{n+2}{n-2}\nabla_{s}\RRR_{ip}\RRR_{jp}-\frac{n+2}{n-2}\RRR_{ip}\nabla_{s}\RRR_{jp}
+\frac{2n}{(n-1)(n-2)}\nabla_{s}\RRR\RRR_{ij}\nonumber\\
&&+\frac{2n}{(n-1)(n-2)}\RRR\nabla_{s}\RRR_{ij}+\frac{2}{n-2}\nabla_{s}\vert\Ric\vert^2g_{ij}-\frac{2}{(n-1)(n-2)}\nabla_{s}\RRR^2
g_{ij}\nonumber\\
&&-2\nabla_{s}\RRR_{kl}\WWW_{kijl}-2\RRR_{kl}\nabla_{s}\WWW_{kijl}+\CCC_{spi}\RRR_{jp}+\CCC_{spj}\RRR_{ip}\nonumber\\
&&+\frac{1}{2(n-1)}\RRR_{jp}[\nabla_i\RRR g_{sp}-\nabla_p\RRR g_{is}]+\frac{1}{2(n-1)}\RRR_{ip}[\nabla_j\RRR g_{sp}-\nabla_p\RRR g_{js}]\,,\nonumber
\end{eqnarray}
where in the last passage we substituted the expression of the Cotton tensor.

We then compute,
\begin{eqnarray*}
\frac{\partial\,}{\partial t}\CCC_{ijk}
&=&\frac{\partial\,}{\partial
  t}\nabla_k\RRR_{ij} -\frac{\partial\,}{\partial t} \nabla_j\RRR_{ik}
-\frac{1}{2(n-1)}\frac{\partial\,}{\partial t}\big(\nabla_k\RRR g_{ij} - \nabla_j\RRR
g_{ik} \big)\\
&=&\nabla_{k}\Delta\RRR_{ij}-\frac{n+2}{n-2}\nabla_{k}\RRR_{ip}\RRR_{jp}-\frac{n+2}{n-2}\RRR_{ip}\nabla_{k}\RRR_{jp}
+\frac{2n}{(n-1)(n-2)}\nabla_{k}\RRR\RRR_{ij}\nonumber\\
&&+\frac{2n}{(n-1)(n-2)}\RRR\nabla_{k}\RRR_{ij}+\frac{2}{n-2}\nabla_{k}\vert\Ric\vert^2g_{ij}-\frac{2}{(n-1)(n-2)}\nabla_{k}\RRR^2
g_{ij}\nonumber\\
&&-2\nabla_{k}\RRR_{pl}\WWW_{pijl}-2\RRR_{pl}\nabla_{k}\WWW_{pijl}+\CCC_{kpi}\RRR_{jp}+\CCC_{kpj}\RRR_{ip}\\
&&+\frac{\RRR_{jp}}{2(n-1)}[\nabla_i\RRR g_{kp}-\nabla_p\RRR g_{ik}]
+\frac{\RRR_{ip}}{2(n-1)}[\nabla_j\RRR g_{kp}-\nabla_p\RRR g_{jk}]\\
&&-\nabla_{j}\Delta\RRR_{ik}+\frac{n+2}{n-2}\nabla_{j}\RRR_{ip}\RRR_{kp}+\frac{n+2}{n-2}\RRR_{ip}\nabla_{j}\RRR_{kp}
-\frac{2n}{(n-1)(n-2)}\nabla_{j}\RRR\RRR_{ik}\nonumber\\
&&-\frac{2n}{(n-1)(n-2)}\RRR\nabla_{j}\RRR_{ik}-\frac{2}{n-2}\nabla_{j}\vert\Ric\vert^2g_{ik}+\frac{2}{n-1}\nabla_{j}\RRR^2
g_{ik}\nonumber\\
&&-2\nabla_{k}\RRR_{pl}\WWW_{pijl}-2\RRR_{pl}\nabla_{k}\WWW_{pijl}-\CCC_{jpi}\RRR_{kp}-\CCC_{jpk}\RRR_{ip}\\
&&-\frac{\RRR_{kp}}{2(n-1)}[\nabla_i\RRR g_{jp}-\nabla_p\RRR g_{ij}]
-\frac{\RRR_{ip}}{2(n-1)}[\nabla_k\RRR g_{jp}-\nabla_p\RRR g_{kj}]\\
&&+\frac{1}{n-1}(\RRR_{ij}\nabla_k\RRR -\RRR_{ik}\nabla_j\RRR\big)\\
&&-\big(\nabla_k\Delta\RRR+2\nabla_{k}|\Ric|^2\big)\frac{g_{ij}}{2(n-1)}
+\big(\nabla_j\Delta\RRR+2\nabla_{j}|\Ric|^2\big) \frac{g_{ik}}{2(n-1)}\\
&=&\nabla_{k}\Delta\RRR_{ij}-\frac{n+2}{n-2}\nabla_{k}\RRR_{ip}\RRR_{jp}-\frac{n+2}{n-2}\RRR_{ip}\nabla_{k}\RRR_{jp}
\\&&+\frac{5n-2}{2(n-1)(n-2)}\nabla_{k}\RRR\RRR_{ij}+\frac{2n}{(n-1)(n-2)}\RRR\nabla_{k}\RRR_{ij}\\
&&+\frac{n}{(n-1)(n-2)}\nabla_{k}\vert\Ric\vert^2g_{ij}-\frac{2}{(n-1)(n-2)}\nabla_{k}\RRR^2
g_{ij}\\
&&+\CCC_{kpi}\RRR_{jp}+\CCC_{kpj}\RRR_{ip}-2\nabla_k \RRR_{pl}\WWW_{pijl}-2\RRR_{pl}\nabla_k\WWW_{pijl}\\
&&-\frac{1}{2(n-1)}\RRR_{pj}\nabla_p\RRR g_{ik}\\
&&-\nabla_{j}\Delta\RRR_{ik}+\frac{n+2}{n-2}\nabla_{j}\RRR_{ip}\RRR_{kp}+\frac{n+2}{n-2}\RRR_{ip}\nabla_{j}\RRR_{kp}
\\&&-\frac{5n-2}{2(n-1)(n-2)}\nabla_{j}\RRR\RRR_{ik}-\frac{2n}{(n-1)(n-2)}\RRR\nabla_{j}\RRR_{ik}\\
&&-\frac{n}{(n-1)(n-2)}\nabla_{j}\vert\Ric\vert^2g_{ik}+\frac{2}{(n-1)(n-2)}\nabla_{j}\RRR^2
g_{ik}\\
&&-\CCC_{jpi}\RRR_{kp}-\CCC_{jpk}\RRR_{ip}+2\nabla_j\RRR_{pl}\WWW_{pikl}+2\RRR_{pl}\nabla_j\WWW_{pikl}\\
&&+\frac{1}{2(n-1)}\nabla_l\RRR\RRR_{lk}g_{ij}-\frac{1}{2(n-1)}\nabla_k\Delta\RRR g_{ij}+\frac{1}{2(n-1)}\nabla_j\Delta\RRR g_{ik}
\end{eqnarray*}
and
$$
\Delta\CCC_{ijk}=\Delta\nabla_k\RRR_{ij}-\Delta\nabla_j\RRR_{ik}-\frac{1}{2(n-1)}\Delta\nabla_k\RRR
g_{ij}+\frac{1}{2(n-1)}\Delta\nabla_j\RRR g_{ik}\,,
$$
hence,
\begin{eqnarray*}
\frac{\partial\,}{\partial t}\CCC_{ijk}-\Delta\CCC_{ijk}
&=&\nabla_{k}\Delta\RRR_{ij}-\nabla_{j}\Delta\RRR_{ik}-\Delta\nabla_k\RRR_{ij}+\Delta\nabla_j\RRR_{ik}\\
&&-\frac{1}{2(n-1)}(\nabla_k\Delta\RRR g_{ij} - \nabla_j\Delta\RRR
g_{ik}-\Delta\nabla_k\RRR g_{ij}+\Delta\nabla_j\RRR g_{ik})\\
&&-\frac{n+2}{n-2}(\nabla_{k}\RRR_{ip}\RRR_{jp}+\RRR_{ip}\nabla_{k}\RRR_{jp})
+\frac{5n-2}{2(n-1)(n-2)}\nabla_{k}\RRR\RRR_{ij}\\
&&+\frac{2n}{(n-1)(n-2)}\RRR\nabla_{k}\RRR_{ij}\\
&&+\frac{n}{(n-1)(n-2)}\nabla_{k}\vert\Ric\vert^2g_{ij}-\frac{2}{(n-1)(n-2)}\nabla_{k}\RRR^2
g_{ij}\\
&&+\CCC_{kpi}\RRR_{jp}+\CCC_{kpj}\RRR_{ip}-2\nabla_k\RRR_{pl}\WWW_{pijl}-2\RRR_{pl}\nabla_k\WWW_{pijl} \\
&&-\frac{1}{2(n-1)}\RRR_{jp}\nabla_p\RRR g_{ik}\\
&&+\frac{n+2}{n-2}(\nabla_{j}\RRR_{ip}\RRR_{kp}+\RRR_{ip}\nabla_{j}\RRR_{kp})
-\frac{5n-2}{2(n-1)(n-2)}\nabla_{j}\RRR\RRR_{ik}\\
&&-\frac{2n}{(n-1)(n-2)}\RRR\nabla_{j}\RRR_{ik}\\
&&-\frac{n}{(n-1)(n-2)}\nabla_{j}\vert\Ric\vert^2g_{ik}+\frac{2}{(n-1)(n-2)}\nabla_{j}\RRR^2
g_{ik}\\
&&-\CCC_{jpi}\RRR_{kp}-\CCC_{jpk}\RRR_{ip}+2\nabla_j\RRR_{pl}\WWW_{pikl}+2\RRR_{pl}\nabla_j\WWW_{pikl}\\
&&+\frac{1}{2(n-1)}\RRR_{kp}\nabla_p\RRR g_{ij}
\end{eqnarray*}

Now to proceed, we need the following commutation rules for the
derivatives of the Ricci tensor and of the scalar curvature, where we
will employ the decomposition formula of the Riemann tensor~\eqref{decriem}.
\begin{eqnarray*}
\nabla_k\Delta\RRR_{ij}-\Delta\nabla_k\RRR_{ij}
&=&\nabla^3_{kll}\RRR_{ij}-\nabla^3_{lkl}\RRR_{ij}+
\nabla^3_{lkl}\RRR_{ij}-\nabla^3_{llk}\RRR_{ij}\\
&=&-\RRR_{kp}\nabla_p\RRR_{ij}
+\RRR_{klip}\nabla_l\RRR_{jp}
+\RRR_{kljp}\nabla_l\RRR_{ip}\\
&&+\nabla^3_{lkl}\RRR_{ij}-\nabla^3_{llk}\RRR_{ij}\\
&=&-\RRR_{kp}\nabla_p\RRR_{ij}
+\frac{1}{2(n-2)}(\RRR_{ik}\nabla_j\RRR
+\RRR_{jk}\nabla_i\RRR)\\
&&-\frac{1}{n-2}(\RRR_{kp}\nabla_i\RRR_{jp}
+\RRR_{kp}\nabla_j\RRR_{ip}
-\RRR_{lp}\nabla_l\RRR_{jp}g_{ik}
-\RRR_{lp}\nabla_l\RRR_{ip}g_{jk})\\
&&-\frac{1}{n-2}(\RRR_{li}\nabla_l\RRR_{jk}
+\RRR_{lj}\nabla_l\RRR_{ik})
-\frac{1}{2(n-1)(n-2)}(\RRR\nabla_j\RRR g_{ik}
+\RRR\nabla_i\RRR g_{jk})\\
&&+\frac{1}{(n-1)(n-2)}(\RRR\nabla_i\RRR_{jk}
+\RRR\nabla_j\RRR_{ik})\\
&&+\nabla_l\big(\RRR_{klip}\RRR_{pj}+\RRR_{kljp}\RRR_{pi}\big)\\
&&+\WWW_{kljp}\nabla_l\RRR_{ip}+\WWW_{klip}\nabla_j\RRR_{jp}\\
&=&-\RRR_{kp}\nabla_p\RRR_{ij}
+\frac{1}{2(n-2)}(\RRR_{ik}\nabla_j\RRR
+\RRR_{jk}\nabla_i\RRR)\\
&&-\frac{1}{n-2}(\RRR_{kp}\nabla_i\RRR_{jp}
+\RRR_{kp}\nabla_j\RRR_{ip}
-\RRR_{lp}\nabla_l\RRR_{jp}g_{ik}
-\RRR_{lp}\nabla_l\RRR_{ip}g_{jk})\\
&&-\frac{1}{n-2}(\RRR_{li}\nabla_l\RRR_{jk}
+\RRR_{lj}\nabla_l\RRR_{ik})
-\frac{1}{2(n-1)(n-2)}(\RRR\nabla_j\RRR g_{ik}
+\RRR\nabla_i\RRR g_{jk})\\
&&+\frac{1}{(n-1)(n-2)}(\RRR\nabla_i\RRR_{jk}
+\RRR\nabla_j\RRR_{ik})\\
&&+\nabla_l\Big(\frac{1}{n-2}(\RRR_{ki}g_{pl}\RRR_{pj}+\RRR_{pl}g_{ki}\RRR_{pj}-\RRR_{li}g_{kp}\RRR_{pj}-\RRR_{kp}g_{li}\RRR_{pj})\\
&&-\frac{1}{(n-1)(n-2)}(\RRR\RRR_{pj}g_{ki}g_{lp}-\RRR\RRR_{pj}g_{kp}g_{il})+\WWW_{klip}\RRR_{pj}\\
&&+\frac{1}{n-2}(\RRR_{kj}g_{lp}\RRR_{pi}+\RRR_{lp}g_{kj}\RRR_{pi}-\RRR_{lj}g_{kp}\RRR_{pi}-\RRR_{kp}g_{lj}\RRR_{pi})\\
&&-\frac{1}{(n-1)(n-2)}(\RRR\RRR_{pi}g_{kj}g_{pl}-\RRR\RRR_{pi}g_{kp}g_{lj})+\WWW_{kljp}\RRR_{pi}\Big)\\
&&+\WWW_{kljp}\nabla_l\RRR_{ip}+\WWW_{klip}\nabla_j\RRR_{jp}\\
&=&-\RRR_{kp}\nabla_p\RRR_{ij}
+\frac{1}{2(n-2)}(\RRR_{ik}\nabla_j\RRR
+\RRR_{jk}\nabla_i\RRR)\\
&&-\frac{1}{n-2}(\RRR_{kp}\nabla_i\RRR_{jp}
+\RRR_{kp}\nabla_j\RRR_{ip}
-\RRR_{lp}\nabla_l\RRR_{jp}g_{ik}
-\RRR_{lp}\nabla_l\RRR_{ip}g_{jk})\\
&&-\frac{1}{n-2}(\RRR_{li}\nabla_l\RRR_{jk}
+\RRR_{lj}\nabla_l\RRR_{ik})
-\frac{1}{2(n-1)(n-2)}(\RRR\nabla_j\RRR g_{ik}
+\RRR\nabla_i\RRR g_{jk})\\
&&+\frac{1}{(n-1)(n-2)}(\RRR\nabla_i\RRR_{jk}
+\RRR\nabla_j\RRR_{ik})\\
&&+\frac{1}{n-2}(\nabla_p\RRR_{ki}\RRR_{pj}+\RRR_{ki}\nabla_j\RRR/2+\nabla_p\RRR g_{ki}\RRR_{pj}/2\\
&&+\RRR_{lp}\nabla_l\RRR_{pj}g_{ik}-\nabla_i\RRR\RRR_{jk}/2-\RRR_{pi}\nabla_p\RRR_{kj}-\nabla_i\RRR_{kp}\RRR_{pj}\\
&&-\RRR_{kp}\nabla_i\RRR_{pj})\\
&&-\frac{1}{(n-1)(n-2)}(\nabla_p\RRR\RRR_{pj}g_{ik}+\RRR\nabla_j\RRR g_{ik}/2-\nabla_i\RRR\RRR_{kj}-\RRR\nabla_i\RRR_{jk})\\
&&+\frac{n-3}{n-2}\CCC_{kip}\RRR_{pj}+\WWW_{klip}\nabla_l\RRR_{pj}\end{eqnarray*}

\begin{eqnarray*}
\phantom{\nabla_k\Delta\RRR_{ij}-\Delta\nabla_k\RRR_{ij}}
&\phantom{=}&+\frac{1}{n-2}(\nabla_p\RRR_{kj}\RRR_{pi}+\RRR_{kj}\nabla_i\RRR/2+\nabla_p\RRR g_{kj}\RRR_{pi}/2\\
&&+\RRR_{lp}g_{kj}\nabla_l\RRR_{pi}-\nabla_j\RRR\RRR_{ki}/2-\RRR_{pj}\nabla_p\RRR_{ki}-\nabla_j\RRR_{kp}\RRR_{pi}-\RRR_{kp}\nabla_j\RRR_{pi})\\
&&-\frac{1}{(n-1)(n-2)}(\nabla_p\RRR\RRR_{pi}g_{kj}+\RRR\nabla_i\RRR g_{jk}/2-\nabla_j\RRR \RRR_{ki}-\RRR\nabla_j\RRR_{ki})\\
&&+\frac{n-3}{n-2}\CCC_{kjp}\RRR_{pi}+\WWW_{kljp}\nabla_l\RRR_{pi}\\
&&+\WWW_{kljp}\nabla_l\RRR_{ip}+\WWW_{klip}\nabla_j\RRR_{jp}\\
&=&-\RRR_{kp}\nabla_p\RRR_{ij}+\frac{n+1}{2(n-1)(n-2)}\RRR_{kj}\nabla_i\RRR-\frac{2}{n-2}\RRR_{kp}\nabla_j\RRR_{ip}\\
&&+\frac{2}{n-2}\RRR_{lp}\nabla_l\RRR_{pi}g_{jk}-\frac{1}{n-2}\RRR_{pj}\nabla_p\RRR_{ik}-\frac{1}{(n-1)(n-2)}\RRR\nabla_i\RRR g_{jk}\\
&&+\frac{2}{(n-1)(n-2)}\RRR\nabla_j\RRR_{ik}+\frac{n+1}{2(n-1)(n-2)}\nabla_j\RRR\RRR_{ki}-\frac{2}{n-2}\RRR_{kp}\nabla_i\RRR_{jp}\\
&&+\frac{2}{n-2}\RRR_{lp}\nabla_l\RRR_{pj}g_{ik}-\frac{1}{n-2}\RRR_{pi}\nabla_p\RRR_{jk}-\frac{1}{(n-1)(n-2)}\RRR\nabla_j\RRR g_{ik}\\
&&+\frac{2}{(n-1)(n-2)}\RRR\nabla_i\RRR_{jk}+2\WWW_{kljp}\nabla_l\RRR_{pi}+2\WWW_{klip}\nabla_l\RRR_{pj}\\
&&+\frac{n-3}{2(n-1)(n-2)}(\nabla_p\RRR g_{ik}\RRR_{pj}+\nabla_p\RRR g_{jk}\RRR_{pi})\\
&&+\frac{n-3}{n-2}(\CCC_{kip}\RRR_{pj}+\CCC_{kjp}\RRR_{pi})-\frac{1}{n-2}(\nabla_i\RRR_{kp}\RRR_{pj}+\nabla_j\RRR_{kp}\RRR_{pi})
\end{eqnarray*}
and 
$$
\nabla_k\Delta\RRR-\Delta\nabla_k\RRR=\RRR_{kllp}\nabla_p\RRR=-\RRR_{kp}\nabla_p\RRR\,.
$$

Then, getting back to the main computation, we obtain
\begin{eqnarray*}
\frac{\partial\,}{\partial t}\CCC_{ijk}-\Delta\CCC_{ijk}
&=&-\RRR_{kp}\nabla_p\RRR_{ij}+\frac{n+1}{2(n-1)(n-2)}\RRR_{kj}\nabla_i\RRR\\
&&-\frac{2}{n-2}\RRR_{kp}\nabla_j\RRR_{ip}+\frac{2}{n-2}\RRR_{lp}\nabla_l\RRR_{pi}g_{jk}\\
&&-\frac{1}{n-2}\RRR_{jp}\nabla_p\RRR_{ik}-\frac{1}{(n-1)(n-2)}\RRR\nabla_i\RRR g_{jk}\\
&&+\frac{2}{(n-1)(n-2)}\RRR\nabla_j\RRR_{ik}+\frac{n+1}{2(n-1)(n-2)}\nabla_j\RRR\RRR_{ki}\\
&&-\frac{2}{n-2}\RRR_{kp}\nabla_i\RRR_{pj}+\frac{2}{n-2}\RRR_{lp}\nabla_l\RRR_{pj}g_{ik}\\
&&-\frac{1}{n-2}\RRR_{pi}\nabla_p\RRR_{kj}-\frac{1}{(n-1)(n-2)}\RRR\nabla_j\RRR g_{ik}\\
&&+\frac{2}{(n-1)(n-2)}\RRR\nabla_i\RRR_{jk}+2\WWW_{kljp}\nabla_l\RRR_{pi}+2\WWW_{klip}\nabla_l\RRR_{pj}\\
&&+\frac{n-3}{2(n-1)(n-2)}(\nabla_p\RRR g_{ik}\RRR_{pj}+\nabla_p\RRR g_{jk}\RRR_{pi})\\
&&+\frac{n-3}{n-2}(\CCC_{kip}\RRR_{pj}+\CCC_{kjp}\RRR_{pi})\\
&&-\frac{1}{n-2}(\nabla_i\RRR_{kp}\RRR_{jp}+\nabla_j\RRR_{kp}\RRR_{pi})\\
&&+\RRR_{jp}\nabla_p\RRR_{ik}-\frac{n+1}{2(n-1)(n-2)}\RRR_{kj}\nabla_i\RRR + \frac{2}{n-2}\RRR_{jp}\nabla_k\RRR_{ip}\\
&&-\frac{2}{n-2}\RRR_{lp}\nabla_l\RRR_{pi} g_{kj}+\frac{1}{n-2}\RRR_{pk}\nabla_p\RRR_{ij}\\
&&+\frac{1}{(n-1)(n-2)}\RRR\nabla_i\RRR g_{jk}-\frac{2}{(n-1)(n-2)}\RRR\nabla_k\RRR_{ij}\\
&&-\frac{n+1}{2(n-1)(n-2)}\nabla_k\RRR\RRR_{ij}+\frac{2}{n-2}\RRR_{jp}\nabla_i\RRR_{kp}\\
&&-\frac{2}{n-2}\RRR_{lp}\nabla_p\RRR_{pk}g_{ij}+\frac{1}{n-2}\RRR_{pi}\nabla_p\RRR_{kj}\\
&&+\frac{1}{(n-1)(n-2)}\RRR\nabla_k\RRR g_{ij}-\frac{2}{(n-1)(n-2)}\RRR\nabla_i\RRR_{kj}\\
&&-2\WWW_{jlkp}\nabla_l\RRR_{pi}-2\WWW_{jlip}\nabla_l\RRR_{pk}\\
&&-\frac{n-3}{2(n-2)(n-2)}(\nabla_p\RRR g_{ij}\RRR_{pk}+\nabla_p\RRR g_{jk}\RRR_{pi})\\
&&-\frac{n-3}{n-2}(\CCC_{jip}\RRR_{pk}+\CCC_{jkp}\RRR_{pi})+\frac{1}{n-2}(\nabla_i\RRR_{pj}\RRR_{pk}+\nabla_k\RRR_{jp}\RRR_{pi})\\
&&+\frac{1}{2(n-1)}(\RRR_{kp}\nabla_p\RRR g_{ij}-\RRR_{jp}\nabla_p\RRR g_{ki})-\frac{n+2}{n-2}(\nabla_k\RRR_{pi}\RRR_{pj}+\RRR_{pi}\nabla_k\RRR_{pj})\\
&&+\frac{n}{(n-1)(n-2)}\nabla_k|\Ric|^2 g_{ij}+\frac{5n-2}{2(n-1)(n-2)}\nabla_k\RRR \RRR_{ij}\\
&&+\frac{2n}{(n-1)(n-2)}\RRR\nabla_k\RRR_{ij}-\frac{2}{(n-1)(n-2)}\nabla_k\RRR^2 g_{ij}\\
&&-2\nabla_k\RRR_{pl}\WWW_{pijl}-2\RRR_{pl}\nabla_k\WWW_{pijl}\\
&&+\CCC_{kli}\RRR_{lj}-\frac{1}{2(n-1)}\nabla_l\RRR\RRR_{lj} g_{ik}+\CCC_{klj}\RRR_{li}\end{eqnarray*}

\begin{eqnarray*}\phantom{\frac{\partial\,}{\partial t}\CCC_{ijk}-\Delta\CCC_{ijk}}&\phantom{=}&+\frac{n+2}{n-2}(\nabla_j\RRR_{pi}\RRR_{pk}+\RRR_{pi}\nabla_j\RRR_{pk})\\
&&-\frac{n}{(n-1)(n-2)}\nabla_j|\Ric|^2g_{ki}-\frac{5n-2}{2(n-1)(n-2)}\nabla_j \RRR\RRR_{ik}\\
&&-\frac{2n}{(n-1)(n-2)}\RRR\nabla_j\RRR_{ik}+\frac{2}{(n-1)(n-2)}\nabla_j\RRR^2 g_{ik}\\
&&+2\nabla_j\RRR_{pl}\WWW_{pikl}+2\RRR_{pl}\nabla_j\WWW_{pikl}\\
&&-\CCC_{jli}\RRR_{lk}+\frac{1}{2(n-1)}\nabla_l\RRR\RRR_{lk}g_{ij}-\CCC_{jlk}\RRR_{li}
\\
&=&\frac{1}{n-2}(\RRR_{pi}\CCC_{jkp}+\RRR_{pk}\CCC_{jip}-\CCC_{kip}\RRR_{pj}-\CCC_{kjp}\RRR_{pi})\\
&&+\Big[\frac{2}{n-2}\RRR_{lp}\nabla_l\RRR_{pj}+\frac{3}{2(n-1)(n-2)}\nabla_j\RRR^2\\
&&-\frac{1}{2(n-2)}\nabla_p\RRR\RRR_{pj}-\frac{n}{(n-1)(n-2)}\nabla_j|\Ric|^2\Big]g_{ik}\\
&&-\Big[\frac{2}{n-2}\RRR_{lp}\nabla_l\RRR_{pk}+\frac{3}{2(n-1)(n-2)}\nabla_k\RRR^2\\
&&-\frac{1}{2(n-2)}\nabla_p\RRR\RRR_{pk}-\frac{n}{(n-1)(n-2)}\nabla_k|\Ric|^2\Big]g_{ij}\\
&&-\frac{n-3}{n-2}\RRR_{kp}\nabla_p\RRR_{ij}+\frac{n-3}{n-2}\RRR_{pj}\nabla_p\RRR_{ik}\\
&&+\frac{n}{n-2}\RRR_{kp}\nabla_j\RRR_{pi}+\frac{n+1}{n-2}\nabla_j\RRR_{pk}\RRR_{pi}-\frac{2}{n-2}\RRR\nabla_j\RRR_{ik}\\
&&-\frac{4n-3}{2(n-1)(n-2)}\nabla_j\RRR\RRR_{ik}-\frac{1}{n-2}\RRR_{kp}\nabla_i\RRR_{pj}+\frac{1}{n-2}\RRR_{pj}\nabla_i\RRR_{pk}\\
&&-\frac{n}{n-2}\RRR_{jp}\nabla_k\RRR_{ip}-\frac{n+1}{n-2}\nabla_k\RRR_{jp}\RRR_{ip}+\frac{2}{n-2}\RRR\nabla_k\RRR_{ij}\\
&&+\frac{4n-3}{2(n-1)(n-2)}\nabla_k\RRR\RRR_{ij}+2\WWW_{klip}\nabla_l\RRR_{pj}+2\WWW_{kljp}\nabla_l\RRR_{pi}-2\WWW_{jlkp}\nabla_l\RRR_{pi}\\
&&-2\WWW_{jlip}\nabla_l\RRR_{pk}-2\nabla_k\RRR_{pl}\WWW_{pijl}-2\RRR_{pl}\nabla_k\WWW_{pijl}+2\nabla_j\RRR_{pl}\WWW_{pikl}+2\RRR_{pl}\nabla_j\WWW_{pikl}
\end{eqnarray*}

Now, by means of the very definition of the Cotton tensor~\eqref{Cottonn}, the 
identities~\eqref{CottonSym}, and the symmetries of the Weyl tensor,
we substitute
\begin{align*}
\CCC_{kpj}-\CCC_{jpk}=&\,-\CCC_{kjp}-\CCC_{jpk}=\CCC_{pkj}\\
\nabla_l\RRR_{jp}=&\,\nabla_j\RRR_{lp}+\CCC_{pjl}+\frac{1}{2(n-1)}
\big(\nabla_l\RRR g_{pj} - \nabla_j\RRR g_{pl} \big)\\
\nabla_l\RRR_{kp}=&\,\nabla_k\RRR_{lp}+\CCC_{pkl}+\frac{1}{2(n-1)}
\big(\nabla_l\RRR g_{pk} - \nabla_k\RRR g_{pl} \big)\\
\nabla_i\RRR_{jp}=&\,\nabla_j\RRR_{ip}+\CCC_{pji}+\frac{1}{2(n-1)}
\big(\nabla_i\RRR g_{jp} - \nabla_j\RRR g_{ip} \big)\\
\nabla_i\RRR_{kp}=&\,\nabla_k\RRR_{ip}+\CCC_{pki}+\frac{1}{2(n-1)}
\big(\nabla_i\RRR g_{kp} - \nabla_k\RRR g_{ip} \big)\\
\nabla_p\RRR_{ij}=&\,\nabla_j\RRR_{pi}+\CCC_{ijp}+\frac{1}{2(n-1)}
\big(\nabla_p\RRR g_{ji} - \nabla_j\RRR g_{pi} \big)\\
\nabla_p\RRR_{ik}=&\,\nabla_k\RRR_{pi}+\CCC_{ikp}+\frac{1}{2(n-1)}
\big(\nabla_p\RRR g_{ki} - \nabla_k\RRR g_{pi} \big)
\end{align*}
in the last expression above, getting
\begin{eqnarray*}
\frac{\partial\,}{\partial t}\CCC_{ijk}-\Delta\CCC_{ijk}
&=&\frac{1}{n-2}(\RRR_{pi}\CCC_{pkj}+\RRR_{pk}\CCC_{jip}-\CCC_{kip}\RRR_{pj})\\
&&+\Big[\frac{2}{n-2}\RRR_{lp}\Big(\nabla_j\RRR_{lp}+\CCC_{pjl}+\frac{1}{2(n-1)}\nabla_l\RRR g_{pj}\\
&&-\frac{1}{2(n-1)}\nabla_j\RRR g_{pl})\Big)+\frac{3}{2(n-1)(n-2)}\nabla_j\RRR^2\\
&&-\frac{1}{2(n-2)}\nabla_p\RRR\RRR_{pj}-\frac{n}{(n-1)(n-2)}\nabla_j|\Ric|^2\Big]g_{ik}\\
&&-\Big[\frac{2}{n-2}\RRR_{lp}\Big(\nabla_k\RRR_{pl}+\CCC_{pkl}+\frac{1}{2(n-1)}\nabla_l\RRR g_{pk}\\
&&-\frac{1}{2(n-1)}\nabla_k\RRR g_{pl}\Big)+\frac{3}{2(n-1)(n-2)}\nabla_k\RRR^2\\
&&-\frac{1}{2(n-2)}\nabla_p\RRR\RRR_{pk}-\frac{n}{(n-1)(n-2)}\nabla_k|\Ric|^2\Big]g_{ij}\\
&&-\frac{n-3}{n-2}\RRR_{kp}\Big(\CCC_{ijp}+\nabla_j\RRR_{ip}+\frac{1}{2(n-1)}(\nabla_p\RRR g_{ij}-\nabla_j\RRR g_{ip})\Big)\\
&&+\frac{n-3}{n-2}\RRR_{pj}\Big(\CCC_{ikp}+\nabla_k\RRR_{ip}+\frac{1}{2(n-1)}(\nabla_p\RRR g_{ik}-\nabla_k\RRR g_{ip})\Big)\\
&&+\frac{n}{n-2}\RRR_{kp}\nabla_j\RRR_{pi}+\frac{n+1}{n-2}\nabla_j\RRR_{pk}\RRR_{pi}-\frac{2}{n-2}\RRR\nabla_j\RRR_{ik}\\
&&-\frac{4n-3}{2(n-1)(n-2)}\nabla_j\RRR\RRR_{ik}\\
&&-\frac{1}{n-2}\RRR_{kp}\Big(\nabla_j\RRR_{ip}+\CCC_{pji}+\frac{1}{2(n-1)}(\nabla_i\RRR g_{jp}-\nabla_j\RRR g_{ip})\Big)\\
&&+\frac{1}{n-2}\RRR_{pj}\Big(\nabla_k\RRR_{ip}+\CCC_{kpi}+\frac{1}{2(n-1)}(\nabla_i\RRR g_{kp}-\nabla_k\RRR g_{ip})\Big)\\
&&-\frac{n}{n-2}\RRR_{jp}\nabla_k\RRR_{ip}-\frac{n+1}{n-2}\nabla_k\RRR_{jp}\RRR_{ip}+\frac{2}{n-2}\RRR\nabla_k\RRR_{ij}\\
&&+\frac{4n-3}{2(n-1)(n-2)}\nabla_k\RRR\RRR_{ij}\\
&&+2\CCC_{plj}\WWW_{pikl}-2\CCC_{plk}\WWW_{pijl}-2\CCC_{pil}\WWW_{jklp}\\
&&-2\WWW_{jklp}\nabla_i\RRR_{pl}-2\RRR_{pl}\nabla_k\WWW_{pijl}+2\RRR_{pl}\nabla_j\WWW_{pikl}\\
&=&\frac{1}{n-2}\left(\RRR_{pi}\CCC_{pkj}+\RRR_{pk}(\CCC_{jip}-\CCC_{pji}-(n-3)\CCC_{ijp})+\RRR_{pj}(\CCC_{pki}-\CCC_{kip}+(n-3)\CCC_{ikp})\right)\\
&&+\frac{2}{n-2}\CCC_{pjl}\RRR_{pl}g_{ik}-\frac{2}{n-2}\CCC_{pkl}\RRR_{pl}g_{ij}-2\CCC_{pjl}\WWW_{pikl}+2\CCC_{pkl}\WWW_{pijl}-2\CCC_{pil}\WWW_{jklp}\\
&&+g_{ik}\Big[\frac{\nabla_j\RRR^2}{(n-1)(n-2)}-\frac{1}{(n-1)(n-2)}\nabla_j|\Ric|^2\Big]\\
&&-g_{ij}\Big[\frac{\nabla_k\RRR^2}{(n-1)(n-2)}-\frac{1}{(n-1)(n-2)}\nabla_k|\Ric|^2\Big]\\
&&-\frac{2}{n-2}\RRR_{jp}\nabla_k\RRR_{ip}-\frac{n+1}{n-2}\nabla_k\RRR_{jp}\RRR_{ip}+\frac{3n-1}{2(n-1)(n-2)}\nabla_k\RRR\RRR_{ij}+\frac{2}{n-2}\RRR\nabla_k\RRR_{ij}\\
&&+\frac{2}{n-2}\RRR_{kp}\nabla_j\RRR_{ip}+\frac{n+1}{n-2}\nabla_j\RRR_{kp}\RRR_{ip}-\frac{3n-1}{2(n-1)(n-2)}\nabla_j\RRR \RRR_{ik}-\frac{2}{n-2}\RRR\nabla_j\RRR_{ik}\\
&&-2\WWW_{jklp}\nabla_i\RRR_{lp}-2\RRR_{lp}\nabla_k\WWW_{pijl}+2\RRR_{pl}\nabla_j\WWW_{pikl}\,.
\end{eqnarray*}
then, we substitute again
\begin{align*}
\nabla_k\RRR_{jp}=&\,\nabla_p\RRR_{kj}+\CCC_{jpk}+\frac{1}{2(n-1)}
\big(\nabla_k\RRR g_{jp} - \nabla_p\RRR g_{jk} \big)\\
\nabla_j\RRR_{kp}=&\,\nabla_p\RRR_{jk}+\CCC_{kpj}+\frac{1}{2(n-1)}
\big(\nabla_j\RRR g_{kp} - \nabla_p\RRR g_{kj} \big)\\
\nabla_k\RRR_{ij}=&\,\nabla_i\RRR_{kj}+\CCC_{jik}+\frac{1}{2(n-1)}
\big(\nabla_k\RRR g_{ij} - \nabla_i\RRR g_{jk} \big)\\
\nabla_j\RRR_{ik}=&\,\nabla_i\RRR_{jk}+\CCC_{kij}+\frac{1}{2(n-1)}
\big(\nabla_j\RRR g_{ik} - \nabla_i\RRR g_{kj} \big)\,,
\end{align*}
finally obtaining
\begin{eqnarray*}
\frac{\partial\,}{\partial t}\CCC_{ijk}-\Delta\CCC_{ijk}
&=&\frac{1}{n-2}\left(\RRR_{pi}\CCC_{pkj}+\RRR_{pk}(\CCC_{jip}-\CCC_{pji}-(n-3)\CCC_{ijp})+\RRR_{pj}(\CCC_{pki}-\CCC_{kip}+(n-3)\CCC_{ikp})\right)\\
&&+\frac{2}{n-2}\CCC_{pjl}\RRR_{pl}g_{ik}-\frac{2}{n-2}\CCC_{pkl}\RRR_{pl}g_{ij}-2\CCC_{pjl}\WWW_{pikl}+2\CCC_{pkl}\WWW_{pijl}-2\CCC_{pil}\WWW_{jklp}\\
&&+g_{ik}\Big[\frac{\nabla_j\RRR^2}{(n-1)(n-2)}-\frac{1}{(n-1)(n-2)}\nabla_j|\Ric|^2\Big]\\
&&-g_{ij}\Big[\frac{\nabla_k\RRR^2}{(n-1)(n-2)}-\frac{1}{(n-1)(n-2)}\nabla_k|\Ric|^2\Big]\\
&&-\frac{2}{n-2}\RRR_{jp}\nabla_k\RRR_{ip}-\frac{n+1}{n-2}\RRR_{ip}\nabla_p\RRR_{kj}-\frac{n+1}{n-2}\RRR_{ip}\CCC_{jpk}\\
&&-\frac{n+1}{2(n-1)(n-2)}\RRR_{ij}\nabla_k\RRR + \frac{n+1}{2(n-1)(n-2)}\RRR_{ip}\nabla_p\RRR g_{jk}\\
&&+\frac{3n-1}{2(n-1)(n-2)}\nabla_k\RRR\RRR_{ij}+\frac{2}{n-2}\RRR(\nabla_i\RRR_{jk}+\CCC_{jik}+\frac{1}{2(n-1)}(\nabla_k\RRR g_{ij}-\nabla_i\RRR g_{jk}))\\
&&+\frac{2}{n-2}\RRR_{kp}\nabla_j\RRR_{ip}+\frac{n+1}{n-2}\RRR_{ip}\nabla_p\RRR_{kj}+\frac{n+1}{n-2}\RRR_{ip}\CCC_{kpj}+\frac{n+1}{2(n-1)(n-2)}\nabla_j\RRR \RRR_{ik}\\
&&-\frac{n+1}{2(n-1)(n-2)}\RRR_{ip}\nabla_p\RRR g_{jk}-\frac{3n-1}{2(n-1)(n-2)}\nabla_j\RRR \RRR_{ik}\\
&&-\frac{2}{n-2}\RRR(\nabla_i\RRR_{jk}+\CCC_{kij}+\frac{1}{2(n-1)}(\nabla_j\RRR g_{ik}-\nabla_i\RRR g_{jk}))\\
&&-2\WWW_{jklp}\nabla_i\RRR_{lp}-2\RRR_{lp}\nabla_k\WWW_{pijl}+2\RRR_{pl}\nabla_j\WWW_{pikl}\\
&=&\frac{1}{n-2}(\RRR_{pk}(\CCC_{jip}-\CCC_{pji}-(n-3)\CCC_{ijp})-\RRR_{pj}(\CCC_{kip}-\CCC_{pki}-(n-3)\CCC_{ikp})\\
&&+(n+2)\RRR_{pi}\CCC_{pkj})+\frac{2}{n-2}(\CCC_{pjl}\RRR_{pl}g_{ik}-\CCC_{pkl}\RRR_{pl}g_{ij})+\frac{2}{n-2}\RRR\CCC_{ijk}\\
&&-2\WWW_{pikl}\CCC_{pjl}+2\WWW_{pijl}\CCC_{pkl}-2\CCC_{pil}\WWW_{jklp}\\
&&+g_{ik}\Big[\frac{\nabla_j\RRR^2}{2(n-1)(n-2)}-\frac{1}{(n-1)(n-2)}\nabla_j|\Ric|^2\Big]\\
&&-g_{ij}\Big[\frac{\nabla_k\RRR^2}{2(n-1)(n-2)}-\frac{1}{(n-1)(n-2)}\nabla_k|\Ric|^2\Big]\\
&&-\frac{2}{n-2}\RRR_{jp}\nabla_k\RRR_{ip}+\frac{1}{n-2}\nabla_k\RRR\RRR_{ij}\\
&&+\frac{2}{n-2}\RRR_{kp}\nabla_j\RRR_{ip}-\frac{1}{n-2}\nabla_j\RRR\RRR_{ik}\\
&&+2\RRR_{lp}\nabla_j\WWW_{pikl}-2\RRR_{lp}\nabla_k\WWW_{pijl}\, ,
\end{eqnarray*}
where in the last passage we used again the 
identities~\eqref{CottonSym} and the fact that 
$$
\WWW_{jklp}\nabla_i\RRR_{lp}=\WWW_{jkpl}\nabla_i\RRR_{pl}=\WWW_{jkpl}\nabla_i\RRR_{lp}=-\WWW_{jklp}\nabla_i\RRR_{lp}\,.
$$

Hence, we can resume this long computation in the following
proposition, getting back to a generic coordinate basis.
\begin{prop}\label{cotn}
During the Ricci flow of a $n$--dimensional Riemannian manifold $(M^n,
g(t))$, 
the Cotton tensor satisfies the following evolution equation
\begin{eqnarray*}
\bigl(\partial_t-\Delta\bigr)\CCC_{ijk}&=&\frac{1}{n-2}\left[g^{pq}\RRR_{pj}(\CCC_{kqi}+\CCC_{qki}+(n-3)\CCC_{ikq})\right.\nonumber\\
&&\left.+(n+2)g^{pq}\RRR_{ip}\CCC_{qkj}
-g^{pq}\RRR_{pk}(\CCC_{jqi}+\CCC_{qji}+(n-3)\CCC_{ijq})\right]\\
&&+\frac{2}{n-2}\RRR\CCC_{ijk}+\frac{2}{n-2}\RRR^{ql}\CCC_{qjl}g_{ik}-\frac{2}{n-2}\RRR^{ql}\CCC_{qkl}g_{ij}\nonumber\\
&&+\frac{1}{(n-1)(n-2)}\nabla_k|\Ric|^2g_{ij}-\frac{1}{(n-1)(n-2)}\nabla_j|\Ric|^2g_{ik}\nonumber\\
&&+\frac{\RRR}{(n-1)(n-2)}\nabla_j\RRR g_{ik}
-\frac{\RRR}{(n-1)(n-2)}\nabla_k\RRR g_{ij}\nonumber\\
&&+\frac{2}{n-2}g^{pq}\RRR_{pk}\nabla_j\RRR_{qi}-\frac{2}{n-2}g^{pq}\RRR_{pj}\nabla_k\RRR_{qi}
+\frac{1}{n-2}\RRR_{ij}\nabla_k\RRR-\frac{1}{n-2}\RRR_{ik}\nabla_j\RRR\,\nonumber\\
&&-2g^{pq}\WWW_{pikl}\CCC_{qjl}+2g^{pq}\WWW_{pijl}\CCC_{qkl}-2g^{pq}\WWW_{jklp}\CCC_{qil}+2g^{pq}\RRR_{pl}\nabla_j\WWW_{qikl}-2g^{pq}\RRR_{pl}\nabla_{k}\WWW_{qijl}\,.\nonumber
\end{eqnarray*}
In particular if the Cotton tensor vanishes identically along the flow we obtain,
\begin{eqnarray*}
0&=&\frac{1}{(n-1)(n-2)}\nabla_k|\Ric|^2g_{ij}-\frac{1}{(n-1)(n-2)}\nabla_j|\Ric|^2g_{ik}\nonumber\\
&&+\frac{\RRR}{(n-1)(n-2)}\nabla_j\RRR g_{ik}
-\frac{\RRR}{(n-1)(n-2)}\nabla_k\RRR g_{ij}\nonumber\\
&&+\frac{2}{n-2}g^{pq}\RRR_{pk}\nabla_j\RRR_{qi}-\frac{2}{n-2}g^{pq}\RRR_{pj}\nabla_k\RRR_{qi}
+\frac{1}{n-2}\RRR_{ij}\nabla_k\RRR-\frac{1}{n-2}\RRR_{ik}\nabla_j\RRR\,\nonumber\\
&&+2g^{pq}\RRR_{pl}\nabla_j\WWW_{qikl}-2g^{pq}\RRR_{pl}\nabla_{k}\WWW_{qijl}\,,\nonumber
\end{eqnarray*}
while, in virtue of relation~\eqref{CottonWeyl}, if the Weyl tensor vanishes along the flow we obtain (compare with~\cite[Proposition~1.1 and Corollary~1.2]{mancat1})
\begin{eqnarray*}
0&=&\frac{1}{(n-1)(n-2)}\nabla_k|\Ric|^2g_{ij}-\frac{1}{(n-1)(n-2)}\nabla_j|\Ric|^2g_{ik}\nonumber\\
&&+\frac{\RRR}{(n-1)(n-2)}\nabla_j\RRR g_{ik}
-\frac{\RRR}{(n-1)(n-2)}\nabla_k\RRR g_{ij}\nonumber\\
&&+\frac{2}{n-2}g^{pq}\RRR_{pk}\nabla_j\RRR_{qi}-\frac{2}{n-2}g^{pq}\RRR_{pj}\nabla_k\RRR_{qi}
+\frac{1}{n-2}\RRR_{ij}\nabla_k\RRR-\frac{1}{n-2}\RRR_{ik}\nabla_j\RRR\,.\nonumber
\end{eqnarray*}
\end{prop}

\begin{cor}\label{CorEvn}
During the Ricci flow of a $n$--dimensional Riemannian manifold $(M^n,
g(t))$, the squared norm of the Cotton tensor satisfies the following evolution equation, in an orthonormal basis,
\begin{eqnarray*}
\bigl(\partial_t-\Delta\bigr)\vert\CCC_{ijk}\vert^2
&=&-2\vert \nabla \CCC_{ijk}\vert^2-\frac{16}{n-2}\CCC_{ipk}\CCC_{iqk}\RRR_{pq}
+\frac{24}{n-2}\CCC_{ipk}\CCC_{kqi}\RRR_{pq}\\
&&+\frac{4}{n-2}\RRR\vert\CCC_{ijk}\vert^2+\frac{8}{n-2}\CCC_{ijk}\RRR_{pk}\nabla_j\RRR_{pi}
+\frac{4}{n-2}\CCC_{ijk}\RRR_{ij}\nabla_k\RRR\nonumber\\
&&+8\CCC_{ijk}\RRR_{lp}\nabla_j\WWW_{pikl}-8\CCC_{ijk}\CCC_{pjl}\WWW_{pikl}-4\CCC_{jpi}\CCC_{ljk}\WWW_{pikl}\,.\nonumber
\end{eqnarray*}
\end{cor}

\begin{eqnarray*}
\bigl(\partial_t-\Delta\bigr)\vert\CCC_{ijk}\vert^2
&=&-2\vert \nabla \CCC_{ijk}\vert^2+2\CCC^{ijk}\RRR_{ip}g^{pq}\CCC_{qjk}+2\CCC^{ijk}\RRR_{jp}g^{pq}\CCC_{iqk}+2\CCC^{ijk}\RRR_{kp}g^{pq}\CCC_{iqk}\nonumber\\
&&+2\CCC_{ijk}\Bigl[\frac{1}{n-2}\left[(\RRR_{pj}(\CCC_{kpi}+\CCC_{pki}+(n-3)\CCC_{ikp})\right.\\
&&+(n+2)\RRR_{pi}\CCC_{pkj}-\RRR_{pk}(\CCC_{jpi}+\CCC_{pji}+(n-3)\CCC_{ijp}))\\
&&+\frac{2}{n-2}\RRR\CCC_{ijk}+\frac{2}{n-2}\RRR_{ql}\CCC_{qjl}g_{ik}-\frac{2}{n-2}\RRR_{ql}\CCC_{qkl}g_{ij}\nonumber\\
&&+\frac{1}{(n-1)(n-2)}\nabla_k|\Ric|^2g_{ij}-\frac{1}{(n-1)(n-2)}\nabla_j|\Ric|^2g_{ik}\nonumber\\
&&+\frac{\RRR}{(n-1)(n-2)}\nabla_j\RRR g_{ik}
-\frac{\RRR}{(n-1)(n-2)}\nabla_k\RRR g_{ij}\nonumber\\
&&+\frac{2}{n-2}\RRR_{qk}\nabla_j\RRR_{qi}-\frac{2}{n-2}\RRR_{qj}\nabla_k\RRR_{qi}
+\frac{1}{n-2}\RRR_{ij}\nabla_k\RRR-\frac{1}{n-2}\RRR_{ik}\nabla_j\RRR\,\nonumber\\
&&-2\WWW_{pikl}\CCC_{pjl}+2\WWW_{pijl}\CCC_{pkl}-2\WWW_{jklp}\CCC_{pil}+2\RRR_{pl}\nabla_j\WWW_{pikl}-2\RRR_{pl}\nabla_{k}\WWW_{pikl}\Bigr]\\
&=&-2\vert \nabla \CCC_{ijk}\vert^2-\frac{16}{n-2}\CCC_{ipk}\CCC_{iqk}\RRR_{pq}
+\frac{24}{n-2}\CCC_{ipk}\CCC_{kqi}\RRR_{pq}\nonumber\\
&&+\frac{4}{n-2}\RRR\vert\CCC_{ijk}\vert^2+\frac{8}{n-2}\CCC_{ijk}\RRR_{pk}\nabla_j\RRR_{pi}
+\frac{4}{n-2}\CCC_{ijk}\RRR_{ij}\nabla_k\RRR\nonumber\\
&&+8\CCC_{ijk}\RRR_{lp}\nabla_j\WWW_{pikl}-8\CCC_{ijk}\CCC_{pjl}\WWW_{pikl}-4\CCC_{jpi}\CCC_{ljk}\WWW_{pikl}\,.\nonumber
\end{eqnarray*}

\begin{rem} Notice that if $n=3$ the two formulas in Proposition~\ref{cotn}
    and Corollary~\ref{CorEvn} become the ones in  Proposition~\ref{cot} and
      Corollary~\ref{CorEv}.
\end{rem}

\section{The Bach Tensor}

The Bach tensor in dimension three is given by
$$
\BBB_{ik}=\nabla_j\CCC_{ijk}\,.
$$

Let $\SSS_{ij}=\RRR_{ij}-\frac{1}{2(n-1)}\Scal g_{ij}$ be the Schouten 
tensor, then 
\begin{equation}\label{eq_DefBach}
\BBB_{ik}=\nabla_j\CCC_{ijk}=\nabla_j(\nabla_k
\SSS_{ij}-\nabla_j\SSS_{ik})=\nabla_j\nabla_k\SSS_{ij}-\Delta\SSS_{ik}\,.
\end{equation}
We compute, in generic dimension $n$,
\begin{eqnarray*}\nabla_{j}\CCC_{ijk}&=&\nabla_j\nabla_k\RRR_{ij}-\frac{1}{2(n-1)}\nabla_j\nabla_k\Scal g_{ij}-\Delta \SSS_{ik}\\
&=&+\RRR_{jkil}\RRR_{jl}+\RRR_{jkjl}\RRR_{il}+\nabla_k\nabla_j\RRR_{ij}-\frac{1}{2(n-1)}\nabla_k\nabla_j\Scal g_{ij}-\Delta\SSS_{ik}\\
&=&+\frac{1}{n-2}\left(\RRR_{ij}g_{kl}-\RRR_{jl}g_{ki}+\RRR_{kl}g_{ij}-\RRR_{ki}g_{jl}-\frac{\Scal}{(n-1)}(g_{ij}g_{kl}-g_{jl}g_{ki})\right)\RRR_{jl}+\WWW_{jkil}\RRR_{jl}\\
&&+\RRR_{kl}\RRR_{il}+\frac{1}{2}\nabla_k\nabla_i\Scal-\frac{1}{2(n-1)}\nabla_k\nabla_i\Scal-\Delta \SSS_{ik}\\
&=&+\frac{1}{n-2}(\RRR_{ji}\RRR_{jk}-|\Ric|^2g_{ik}+\RRR_{kl}\RRR_{il}-\Scal\RRR_{ik})-\frac{\Scal}{(n-1)(n-2)}\RRR_{ik}+\frac{\Scal^2}{(n-1)(n-2)}g_{ik}\\
&&+\WWW_{jkil}\RRR_{jl}+\RRR_{kl}\RRR_{il}+\frac{n-2}{2(n-1)}\nabla_k\nabla_i\Scal-\Delta\SSS_{ik}\\
&=&\frac{n}{n-2}\RRR_{ij}\RRR_{kj}-\frac{n}{(n-1)(n-2)}\Scal\RRR_{ik}-\frac{1}{n-2}|\Ric|^2g_{ik}+\frac{\Scal^2}{(n-1)(n-2)}g_{ik}\\&&+\WWW_{jkil}\RRR_{jl}+\frac{n-2}{2(n-1)}\nabla_k\nabla_i\Scal-\Delta \SSS_{ik}\,.
\end{eqnarray*}
From this last expression, it is easy to see that the Bach tensor in dimension $3$ is
symmetric, i.e. $\BBB_{ik}=\BBB_{ki}$. Moreover, it is trace--free, that is, $g^{ik}\BBB_{ik}=0$ as $g^{ik}\nabla\CCC_{ijk}=0$.

\begin{rem}In higher dimension, the Bach tensor is given by
$$\BBB_{ik}=\frac{1}{n-2}(\nabla_j \CCC_{ijk}-\RRR_{jl}\WWW_{ijkl})\,.$$
We note that, since
$\RRR_{jl}\WWW_{ijkl}=\RRR_{jl}\WWW_{klij}=\RRR_{jl}\WWW_{kjil}$, from the above computation we get that the Bach tensor is symmetric in any dimension; finally, as the
Weyl tensor is trace-free in every pair of indexes, there holds
$g^{ik}\BBB_{ik}=0$.
\end{rem}

We recall that Schur lemma yields the following equation for the divergence of the Schouten tensor
\begin{equation}\label{eq_SchoutenSchur}\nabla_j\SSS_{ij}=\frac{n-2}{2(n-1)}\nabla_i\Scal\,.\end{equation}
We write
$$
\nabla_k\nabla_j\CCC_{ijk}=\nabla_k\nabla_j\nabla_k\SSS_{ij}-\nabla_k\nabla_j\nabla_j\SSS_{ik}=[\nabla_j,\nabla_k]\nabla_j\SSS_{ik}\,,
$$
therefore,
\begin{eqnarray*}\nabla_k\nabla_j\CCC_{ijk}&=&\RRR_{jkjl}\nabla_l\SSS_{ik}+\RRR_{jkil}\nabla_j\SSS_{lk}+\RRR_{jkkl}\nabla_j\SSS_{li}\\
&=&\RRR_{kl}\nabla_l\SSS_{ik}+\RRR_{jkil}\nabla_j\SSS_{lk}-\RRR_{jl}\nabla_j\SSS_{li}\\
&=&\left[\frac{1}{n-2}\,(\RRR_{ij}g_{kl}-\RRR_{jl}g_{ik}+\RRR_{kl}g_{ij}-\RRR_{ik}g_{jl})-\frac{1}{(n-1)(n-2)}\,\RRR(g_{ij}g_{kl}-g_{ik}g_{jl}) + \WWW_{jkil}\right]\nabla_j\SSS_{lk}\\
&=& \frac{1}{n-2}(-\RRR_{jl}\nabla_j\SSS_{il}+\RRR_{kl}\nabla_i\SSS_{kl})+\WWW_{jkil}\nabla_j\SSS_{lk}\\
&=&\frac{1}{n-2}\RRR_{jl}(\nabla_i\SSS_{lj}-\nabla_j\SSS_{il})+\WWW_{jkil}\nabla_j\RRR_{kl}\\
&=&\frac{1}{n-2}\RRR_{jl}\CCC_{lji}+\WWW_{iljk}\nabla_j\RRR_{kl}\,,
\end{eqnarray*}
where we repeatedly used equation~\eqref{eq_SchoutenSchur}, the trace--free property  of the Weyl tensor and the definition of the Cotton tensor.

Recalling that
$$
\nabla_k\WWW_{ijkl}=\nabla_k\WWW_{klij}=-\frac{n-3}{n-2}\CCC_{lij}=\frac{n-3}{n-2}\CCC_{lji}\,,
$$
the divergence of the Bach tensor is given by
\begin{eqnarray*}\nabla_k\BBB_{ik}&=&\frac{1}{n-2}\nabla_k(\nabla_j \CCC_{ijk}-\RRR_{jl}\WWW_{ijkl})=\frac{1}{(n-2)^2}\RRR_{jl}\CCC_{jli}-\frac{n-3}{(n-2)^2}\,\CCC_{jli}\RRR_{jl}\\&=&-\frac{n-4}{(n-2)^2}\CCC_{jli}\RRR_{jl}\,.\end{eqnarray*}
In particular, for $n=3$, we obtain $\nabla_k\BBB_{ik}=\nabla_k\BBB_{ki}=\RRR_{jl}\CCC_{jli}$ and, for $n=4$, we get the classical result $\nabla_k\BBB_{ik}=\nabla_k\BBB_{ki}=0$.

\subsection{The Evolution Equation of the Bach Tensor in 3D}\ \\

We turn now our attention to the evolution of the Bach tensor along
the Ricci flow in dimension three. In order to obtain its evolution
equation, instead of calculating directly the time derivative and the
Laplacian of the Bach tensor, we employ the following equation
\begin{equation}\label{eq_EvBach}
(\partial_t-\Delta)\BBB_{ik}=\nabla_j(\partial_t-\Delta)\CCC_{ijk}-[\Delta,\nabla_j]\CCC_{ijk}+2\RRR_{pj}\nabla_p\CCC_{ijk}+[
\partial_t, \nabla_j]C_{ijk}\,,
\end{equation}
which relates the quantity we want to compute with the evolution of the
Cotton tensor, the evolution of the Christoffel symbols and the
formulas for the exchange of covariant derivatives. 
We will work on the various terms separately.

By the commutations formulas for derivatives, we have

\begin{eqnarray*}
\nabla_l\nabla_l\nabla_q \CCC_{ijk}-\nabla_l\nabla_q\nabla_l\CCC_{ijk}=\nabla_l(\RRR_{lqip}\CCC_{pjk}+\RRR_{lqjp}\CCC_{ipk}+\RRR_{lqkp}\CCC_{ijp})\end{eqnarray*}
\begin{eqnarray*}\nabla_l\nabla_q\nabla_s\CCC_{ijk}-\nabla_q\nabla_l\nabla_s\CCC_{ijk}=\RRR_{lqsp}\nabla_p\CCC_{ijk}+\RRR_{lqip}\nabla_s\CCC_{pjk}+\RRR_{lqjp}\nabla_s\CCC_{ipk}+\RRR_{lqkp}\nabla_s\CCC_{ijp},
\end{eqnarray*}
and putting these together with $q=j$ and $l=s$, we get
\begin{eqnarray*}
[\Delta,\nabla_j]\CCC_{ijk}&=&\nabla_l(\RRR_{ljip}\CCC_{pjk}-\RRR_{lp}\CCC_{ipk}+\RRR_{ljkp}\CCC_{ijp})\\
&&+\RRR_{jp}\nabla_p\CCC_{ijk}+\RRR_{ljip}\nabla_l\CCC_{pjk}-\RRR_{lp}\nabla_l\CCC_{ipk}+\RRR_{ljkp}\nabla_l\CCC_{ijp}\\
&=&\nabla_l\left[\left(\RRR_{li}g_{jp}-\RRR_{lp}g_{ji}+\RRR_{jp}g_{li}-\RRR_{ji}g_{lp}-\frac{\RRR}{2}(g_{li}g_{jp}-g_{lp}g_{ji})\right)\CCC_{pjk}\right.\\
&&\left.-\RRR_{lp}\CCC_{ipk}+\left(\RRR_{lk}g_{jp}-\RRR_{lp}g_{jk}+\RRR_{jp}g_{lk}-\RRR_{jk}g_{lp}-\frac{\RRR}{2}(g_{lk}g_{jp}-g_{lp}g_{jk})\right)\CCC_{ijp}\right]\\
&&+\RRR_{jp}\nabla_{p}\CCC_{ijk}+\RRR_{ljip}\nabla_{l}\CCC_{pjk}-\RRR_{lp}\nabla_{l}\CCC_{ipk}+\RRR_{ljkp}\nabla_{l}\CCC_{ijp}\\
&=&-\frac{1}{2}\nabla_{p}\RRR\CCC_{pik}-\RRR_{lp}\nabla_{l}\CCC_{pik}+\nabla_{i}\RRR_{jp}\CCC_{pjk}+\RRR_{jp}\nabla_{i}\CCC_{pjk}-\nabla_{p}\RRR_{ji}\CCC_{pjk}-\RRR_{ji}\nabla_{p}\CCC_{pjk}\\
&&+\frac{1}{2}\nabla_{p}\RRR\CCC_{pik}+\frac{\RRR}{2}\nabla_{p}\CCC_{pik}-\frac{1}{2}\nabla_{p}\RRR\CCC_{ipk}-\RRR_{lp}\nabla_{l}\CCC_{ipk}-\frac{1}{2}\nabla_{p}\RRR\CCC_{ikp}-\RRR_{lp}\nabla_{l}\CCC_{ikp}\\
&&+\nabla_{k}\RRR_{jp}\CCC_{ijp}+\RRR_{jp}\nabla_{k}\CCC_{ijp}-\nabla_{p}\RRR_{jk}\CCC_{ijp}-\RRR_{jk}\nabla_{p}\CCC_{ijp}+\frac{1}{2}\nabla_{p}\RRR\CCC_{ikp}+\frac{\RRR}{2}\nabla_{p}\CCC_{ikp}\\
&&+\RRR_{jp}\nabla_{p}\CCC_{ijk}-\RRR_{lp}\nabla_{l}\CCC_{pik}+\RRR_{jp}\nabla_{i}\CCC_{pjk}-\RRR_{ji}\nabla_{p}\CCC_{pjk}+\frac{\RRR}{2}\nabla_{p}\CCC_{pik}\\
&&-\RRR_{lp}\nabla_{l}\CCC_{ipk}-\RRR_{lp}\nabla_{l}\CCC_{ikp}+\RRR_{jp}\nabla_{k}\CCC_{ijp}-\RRR_{jk}\nabla_{p}\CCC_{ijp}+\frac{\RRR}{2}\nabla_{p}\CCC_{ikp}\\
&=&\nabla_{i}\RRR_{jp}\CCC_{pjk}-\nabla_{p}\RRR_{ji}\CCC_{ijp}-\nabla_{p}\RRR_{jk}\CCC_{ijp}-\nabla_{p}\RRR_{jk}\CCC_{ijp}-2\RRR_{lp}\nabla_{l}\CCC_{pik}\\
&&+2\RRR_{lp}\nabla_{i}\CCC_{plk}-2\RRR_{ji}\nabla_{p}\CCC_{pjk}+\RRR\nabla_{p}\CCC_{pik}+\frac{1}{2}\nabla_{p}\RRR\CCC_{ikp}+2\RRR_{jp}\nabla_{k}\CCC_{ijp}\\
&&-2\RRR_{jk}\nabla_{p}\CCC_{ijp}+\RRR\nabla_{p}\CCC_{ikp}+\RRR_{jp}\nabla_{p}\CCC_{ijk}\\
&=&\nabla_{i}\RRR_{lp}\CCC_{plk}-\nabla_{p}\RRR_{li}\CCC_{plk}+\nabla_{k}\RRR_{lp}\CCC_{ilp}-\nabla_{p}\RRR_{lk}\CCC_{ilp}\\
&&-2\RRR_{lp}\nabla_{l}\CCC_{pik}+2\RRR_{lp}\nabla_{i}\CCC_{plk}+2\RRR_{li}\BBB_{kl}-2\RRR_{li}\BBB_{lk}+2\RRR_{lp}\nabla_{k}\CCC_{ilp}\\
&&+2\RRR_{lk}\BBB_{il}+\RRR_{lp}\nabla_{p}\CCC_{ilk}-\RRR\BBB_{ik}+\frac{1}{2}\nabla_{p}\RRR\CCC_{ikp}+\RRR\BBB_{ik}-\RRR\BBB_{ik}\\
&=&\nabla_{i}\RRR_{lp}\CCC_{plk}-\nabla_{p}\RRR_{li}\CCC_{plk}+\nabla_{k}\RRR_{lp}\CCC_{ilp}-\nabla_{p}\RRR_{lk}\CCC_{ilp}\\
&&+\RRR_{lp}\nabla_{p}\CCC_{ilk}+2\RRR_{lp}\nabla_{i}\CCC_{plk}+2\RRR_{lp}\nabla_{k}\CCC_{ilp}-2\RRR_{lp}\nabla_{l}\CCC_{ipk}\\
&&+\frac{1}{2}\nabla_{p}\RRR\CCC_{ikp}+2\RRR_{lk}\BBB_{il}-\RRR\BBB_{ik}\,.
\end{eqnarray*}
The covariant derivative of the evolution of the Cotton tensor is given by
\begin{eqnarray*}\nabla_j(\partial_t-\Delta)\CCC_{ijk}&=&\frac{5}{2}\nabla_{p}\RRR\CCC_{ipk}+\nabla_{p}\RRR\CCC_{pki}+\RRR_{lp}\nabla_{p}\CCC_{kli}+\RRR_{lp}\nabla_{p}\CCC_{lki}-\nabla_{p}\RRR_{kl}\CCC_{pli}\\
&&-\nabla_{p}\RRR_{kl}\CCC_{lpi}-\RRR_{kp}\BBB_{pi}+5\nabla_{p}\RRR_{il}\CCC_{lkp}-5\RRR_{ip}\BBB_{pk}+2\RRR\BBB_{ik}\\
&&+2\nabla_{s}\RRR_{pl}\CCC_{psl}g_{ik}+2\RRR_{pl}\BBB_{pl}g_{ik}-2\nabla_{i}\RRR_{pl}\CCC_{pkl}-2\RRR_{pl}\nabla_{i}\CCC_{pkl}\\
&&+\frac{1}{2}(|\nabla\RRR|^2+\RRR\Delta\RRR-\Delta|\Ric|^2)g_{ik}-\frac{1}{2}(\nabla_{i}\RRR\nabla_{k}\RRR+\RRR\nabla_{i}\nabla_{k}\RRR-\nabla_{i}\nabla_{k}|\Ric|^2)\\
&&+2\Delta\RRR_{ip}\RRR_{kp}+2\nabla_{l}\RRR_{ip}\nabla_{l}\RRR_{kp}-2\nabla_{l}\nabla_{k}\RRR_{ip}\RRR_{lp}-\nabla_{k}\RRR_{ip}\nabla_{p}\RRR\\
&&+\nabla_{l}\nabla_{k}\RRR\RRR_{il}+\frac{1}{2}\nabla_{k}\RRR\nabla_{i}\RRR-\Delta\RRR\RRR_{ik}-\nabla_{l}\RRR\nabla_{l}\RRR_{ik}.
\end{eqnarray*}
Finally, the commutator between the covariant derivative and the time derivative can be expressed in terms of the time derivatives of the Christoffel symbols, as follows
\begin{eqnarray*}
[\partial_t, \nabla_j]\CCC_{ijk}&=&-\partial_t\Gamma_{ij}^p\CCC_{pjk}-\partial_t\Gamma_{jk}^p\CCC_{ijp}\\
&=&\nabla_i\RRR_{jp}\CCC_{pjk}+\nabla_j\RRR_{ip}\CCC_{pjk}-\nabla_p\RRR_{ij}\CCC_{pjk}+\nabla_j\RRR_{kp}\CCC_{ijp}+\nabla_k\RRR_{jp}\CCC_{ijp}-\nabla_{p}\RRR_{jk}\CCC_{ijp}\\
&=&\nabla_{i}\RRR_{jp}\CCC_{pjk}+\nabla_{p}\RRR_{ij}\CCC_{jpk}+\nabla_{p}\RRR_{ij}\CCC_{pkj}+\nabla_{p}\RRR_{kj}\CCC_{ipj}+\nabla_{k}\RRR_{jp}\CCC_{ijp}+\nabla_{p}\RRR_{jk}\CCC_{ipj}\\
&=&\nabla_{i}\RRR_{jp}\CCC_{pjk}-\nabla_{p}\RRR_{ij}\CCC_{pkj}-\nabla_{p}\RRR_{ij}\CCC_{kjp}+\nabla_{p}\RRR_{ij}\CCC_{pkj}+2\nabla_{p}\RRR_{kj}\CCC_{ipj}\\
&=&\nabla_{i}\RRR_{jp}\CCC_{pjk}-\nabla_{p}\RRR_{ij}\CCC_{kjp}+2\nabla_{p}\RRR_{kj}\CCC_{ipj}\,.
\end{eqnarray*}
Substituting into~\eqref{eq_EvBach}, and making some computations, we obtain the evolution equation
\begin{prop}\label{EvBachRF3D}
During the Ricci flow of a $3$--dimensional Riemannian manifold $(M^3, g(t))$ the Bach tensor satisfies the following evolution equation
\begin{eqnarray*}
(\partial_t-\Delta)\BBB_{ik}=&&\left[3\nabla_{p}\RRR\CCC_{ipk}+\nabla_{p}\RRR\CCC_{pki}-\nabla_{p}\RRR\nabla_{k}\RRR_{ip}\right]\\
&+&\left[-2\RRR_{pl}\nabla_{p}\CCC_{ikl}-3\RRR_{pk}\BBB_{pi}-5\RRR_{pi}\BBB_{pk}+2\Delta\RRR_{ip}\RRR_{kp}\right.\\
&&\,\left.-2\nabla_{l}\nabla_{k}\RRR_{pi}\RRR_{pl}+\nabla_{l}\nabla_{k}\RRR\RRR_{li}-\Delta\RRR\RRR_{ik}\right]\\
&+&\left[-2\nabla_{p}\RRR_{kl}\CCC_{lpi}-2\nabla_{p}\RRR_{kl}\CCC_{ilp}-4\nabla_{p}\RRR_{il}\CCC_{lpk}-2\nabla_{i}\RRR_{pl}\CCC_{pkl}\right]\\
&+&\left[3\RRR\BBB_{ik}+2\nabla_{s}\RRR_{pl}\CCC_{psl}g_{ik}+2\RRR_{pl}\BBB_{pl}g_{ik}\right.\\
&&\,\left.+\frac{1}{2}(|\nabla\RRR|^2+\RRR\Delta\RRR-\Delta|\Ric|^2)g_{ik}-\frac{1}{2}(\RRR\nabla_{i}\nabla_{k}\RRR-\nabla_{i}\nabla_{k}|\Ric|^2)\right.\\
&&\,\left.+2\nabla_{l}\RRR_{ip}\nabla_{l}\RRR_{kp}-\nabla_{l}\RRR\nabla_{l}\RRR_{ik}\right].
\end{eqnarray*}

Hence, if the Bach tensor vanishes identically along the flow, we have
\begin{eqnarray*}
0&=&3\nabla_{p}\RRR\CCC_{ipk}+\nabla_{p}\RRR\CCC_{pki}-\nabla_{p}\RRR\nabla_{k}\RRR_{ip}-2\RRR_{pl}\nabla_{p}\CCC_{ikl}\\
&&+2\Delta\RRR_{ip}\RRR_{kp}-2\nabla_{l}\nabla_{k}\RRR_{pi}\RRR_{pl}+\nabla_{l}\nabla_{k}\RRR\RRR_{li}-\Delta\RRR\RRR_{ik}\\
&&-2\nabla_{p}\RRR_{kl}\CCC_{lpi}-2\nabla_{p}\RRR_{kl}\CCC_{ilp}-4\nabla_{p}\RRR_{il}\CCC_{lpk}-2\nabla_{i}\RRR_{pl}\CCC_{pkl}\\
&&+2\nabla_{s}\RRR_{pl}\CCC_{psl}g_{ik}+\frac{1}{2}(|\nabla\RRR|^2+\RRR\Delta\RRR-\Delta|\Ric|^2)g_{ik}\\
&&-\frac{1}{2}(\RRR\nabla_{i}\nabla_{k}\RRR-\nabla_{i}\nabla_{k}|\Ric|^2)+2\nabla_{l}\RRR_{ip}\nabla_{l}\RRR_{kp}-\nabla_{l}\RRR\nabla_{l}\RRR_{ik}.
\end{eqnarray*}
\end{prop}

\begin{rem}
Note that, from the symmetry property of the Bach tensor, we have that
the RHS in the evolution equation of the Bach tensor should be
symmetric in the two indices. It is not so difficult to check that
this property is verified for the formula 
in Proposition \ref{EvBachRF3D}. Indeed, each of the terms in between
square brackets is symmetric in the two indices.
\end{rem}

As a consequence of Proposition \ref{EvBachRF3D}, we get that during the Ricci flow  of a $3$--dimensional Riemannian manifold the squared norm of the Bach tensor satisfies
\begin{eqnarray*}(\partial_t-\Delta)|\BBB_{ik}|^2&=&-2|\nabla \BBB_{ik}|^2-12\BBB_{ik}\BBB_{iq}\RRR_{qk}+6\BBB_{ik}\nabla_{p}\RRR-4\BBB_{ik}\RRR_{pl}\nabla_{p}\CCC_{ikl}\\
&&+4\BBB_{ik}\nabla_{p}\RRR_{kl}\CCC_{pil}-8\BBB_{ik}\nabla_{p}\RRR_{kl}\CCC_{lpi}-4\BBB_{ik}\nabla_{i}\RRR_{pl}\CCC_{pkl}+6\RRR|\BBB_{ik}|^2\\
&&-2\BBB_{ik}\nabla_{p}\RRR\nabla_{k}\RRR_{ip}+4\BBB_{ik}\Delta\RRR_{ip}\RRR_{kp}-4\BBB_{ik}\nabla_{l}\nabla_{k}\RRR_{pi}\RRR_{pl}+2\BBB_{ik}\nabla_{l}\nabla_{k}\RRR\RRR_{li}\\
&&-2\BBB_{ik}\Delta\RRR\RRR_{ik}-\BBB_{ik}\RRR\nabla_{i}\nabla_{k}\RRR+\BBB_{ik}\nabla_{i}\nabla_{k}|\Ric|^2-2\BBB_{ik}\nabla_{l}\RRR\nabla_{l}\RRR_{ik}\\
&&+4\BBB_{ik}\nabla_{l}\RRR_{ip}\nabla_{l}\RRR_{kp}.
\end{eqnarray*}

\subsection{The Bach Tensor of Three--Dimensional Gradient Ricci Solitons}\ \\

In what follows, we will use formulas~\eqref{SolEq0}--\eqref{SolEq4} to
derive an expression of the Bach tensor and of its divergence in the
particular case of a gradient Ricci soliton in dimension three.

By straightforward computations, we obtain
\begin{eqnarray*}\BBB_{ik}&=&\nabla_j \CCC_{ijk}\\
&=&\frac{\nabla_i\nabla_k \Scal}{4}-\frac{\Delta\Scal}{4}g_{ik}-\nabla_j\RRR_{ik}\nabla_j f + \frac{g_{ik}}2 \nabla_j\Scal\nabla_j f\\
&&+\left(\RRR_{ij}-\frac{\Scal}{2}g_{ij}\right)\nabla_j\nabla_k f - \left(\RRR_{ik}-\frac{\Scal}{2}g_{ik}\right) \Delta f\\
&=&\frac{1}{4}\nabla_i\nabla_k \Scal - \frac{1}{4}\Delta \Scal g_{ik}-\nabla_j\RRR_{ik}\nabla_j f + \frac{1}{2}\nabla_j\Scal\nabla_jf g_{ik}-\RRR_{ij}\RRR_{jk}+\lambda \RRR_{ik}\\
&&+\frac{1}{2}\Scal\RRR_{ik}-\frac{\lambda}{2}\Scal g_{ik}-3\lambda\RRR_{ik}+\Scal\RRR_{ik}+\frac32\lambda\Scal g_{ik}-\frac{1}2\Scal^2g_{ik}\\
&=&\frac{1}{2}\nabla_i\RRR_{lk}\nabla_l f-\frac{1}{2}\RRR_{lk}\RRR_{li}+\frac{\lambda}{2}\RRR_{ik}-\frac{1}{4}\nabla_l\Scal\nabla_l fg_{ik}-\frac{\lambda}{2}\Scal g_{ik}\\
&&+\frac{1}{2}|\Ric|^2g_{ik}-\nabla_j\RRR_{ik}\nabla_j f + \frac{1}{2}\nabla_j\Scal\nabla_j f g_{ik}-\RRR_{ij}\RRR_{jk}+\lambda\RRR_{ik}+\frac{1}{2}\Scal\RRR_{ik}\\
&&-\frac{\lambda}{2}\Scal g_{ik}-3\lambda\RRR_{ik}+\Scal\RRR_{ik}+\frac{3}{2}\lambda\Scal g_{ik}-\frac{1}{2}\Scal^2 g_{ik}\\
&=&\frac{1}{2}\nabla_i\RRR_{lk}\nabla_l f + \frac{1}{4}\nabla_j\Scal\nabla_j f g_{ik}-\nabla_j\RRR_{ik}\nabla_j f -\frac{3}{2}\RRR_{ij}\RRR_{jk}-\frac{3}{2}\lambda\RRR_{ik}\\
&&+\frac{3}{2}\Scal\RRR_{ik}+\frac{\lambda}{2}\Scal g_{ik}+\frac{1}{2}|\Ric|^2 g_{ik}-\frac{1}{2}\Scal^2 g_{ik}\,.
\end{eqnarray*}
A more compact formulation, employing equations~\eqref{SolEq1} and~\eqref{SolEq2}, is given by
\begin{equation*}
\BBB_{ik}=\frac{1}{2}\nabla_i\RRR_{lk}\nabla_l f+\frac{1}{4}\Delta
\Scal g_{ik}-\frac{1}{2}\Delta
\RRR_{ik}-\frac{1}{2}\nabla_j\RRR_{ik}\nabla_j
f-\frac{\lambda}{2}\RRR_{ik}+\frac{1}{2}\RRR_{ij}\RRR_{jk}\,.
\end{equation*}
Moreover, as we know that $\nabla_k\BBB_{ik}=\CCC_{lji}\RRR_{lj}$, we have
\begin{eqnarray*}
\nabla_{k}\BBB_{ik}&=&\frac{1}{4}\Scal \nabla_i\Scal-\frac{1}{4}\RRR_{ij}\nabla_j\Scal + |\Ric|^2\nabla_i f -\frac{1}{2}\Scal^2\nabla_i f-\RRR_{il}\nabla_j f\RRR_{lj}+\frac{1}{2}\Scal\RRR_{ij}\nabla_j f\\
&=&\frac{1}{2}\Scal\nabla_i\Scal-\frac{3}{4}\RRR_{il}\nabla_l\Scal +
|\Ric|^2\nabla_i f- \frac{1}{2}\Scal^2\nabla_i f\,.
\end{eqnarray*}
Therefore, if the divergence of the Bach tensor vanishes, we conclude
\begin{equation*}
\frac{1}{2}\Scal\nabla_i\Scal-\frac{3}{4}\RRR_{ik}\nabla_k\Scal +
|\Ric|^2\nabla_i f- \frac{1}{2}\Scal^2\nabla_i f=0\,.
\end{equation*}
Taking the scalar product with $\nabla f$ in both sides of this equation, we obtain
$$
0=\frac{1}{2}\Scal\langle\nabla\Scal, \nabla f\rangle -
\frac{3}{8}|\nabla \Scal|^2 + |\Ric|^2|\nabla f|^2 -
\frac{1}{2}\Scal^2|\nabla f|^2
$$
and, from formulas~\eqref{SolEq4} and~\eqref{eq_NormaCottonSol}, we compute
\begin{eqnarray*}|\CCC_{ijk}|^2&=&(\RRR_{ij}\nabla_k f-\RRR_{ik}\nabla_j f)\left(\frac{\nabla_k \Scal}{4}g_{ij}-\frac{\nabla_j \Scal}{4}g_{ik}+\left(\RRR_{ij}-\frac{\Scal}{2}g_{ij}\right)\nabla_k f - \left(\RRR_{ik}-\frac{\Scal}{2}g_{ik}\right)\nabla_j f\right)\\
&=&\frac{\Scal}{4}\nabla_k\Scal\nabla_k f-\frac{1}{4}\RRR_{kj}\nabla_j \Scal\nabla_k f+|\Ric|^2|\nabla f|^2-\frac{\Scal^2}{2}|\nabla f|^2-\RRR_{ij}\nabla_j f\RRR_{ik}\nabla_k f+\frac{\Scal}{2}\RRR_{kj}\nabla_k f\nabla_j f\\ & &- \frac{1}{4}\RRR_{jk}\nabla_j f\nabla_k\Scal+\frac{\Scal}{4}\nabla_j\Scal\nabla_j f-\RRR_{ik}\nabla_k f\RRR_{ij}\nabla_j f+\frac{\Scal}{2}\RRR_{jk}\nabla_j f\nabla_k f+|\Ric|^2|\nabla f|^2 - \frac{\Scal^2}{2}|\nabla f|^2\\
&=&2|\Ric|^2|\nabla f|^2 - \Scal^2|\nabla f|^2+\Scal\nabla_k\Scal\nabla_k f-\frac{3}{4}|\nabla \Scal|^2\,,
\end{eqnarray*}
where we repeatedly used equation~\eqref{SolEq3}. \\
Therefore, we obtain 
$$
\nabla_k \BBB_{ik}\nabla_i f=\frac{1}{2}|\CCC_{ijk}|^2\,,
$$
so, if the divergence of the Bach tensor vanishes then the
Cotton tensor vanishes as well (this was already obtained in~\cite{caochen3}). 
As a consequence, getting back to Section~\ref{cgrad}, the soliton 
is locally a warped product of a constant curvature surface on a
interval of $\RR$.

\bibliographystyle{amsplain}
\bibliography{biblio}

\providecommand{\bysame}{\leavevmode\hbox to3em{\hrulefill}\thinspace}
\providecommand{\MR}{\relax\ifhmode\unskip\space\fi MR }
\providecommand{\MRhref}[2]{%
  \href{http://www.ams.org/mathscinet-getitem?mr=#1}{#2}
}
\providecommand{\href}[2]{#2}
\begin{thebibliography}{1}

\bibitem{besse}
A.~L. Besse, \emph{Einstein manifolds}, Springer--Verlag, Berlin, 2008.

\bibitem{caochen3}
H.-D. Cao, G.~Catino, Q.~Chen, C.~Mantegazza, and L.~Mazzieri, \emph{Bach--flat
  gradient steady {R}icci solitons}, Calc. Var. Partial Differential Equations
  \textbf{49} (2014), no.~1-2, 125--138.

\bibitem{caochenzhu}
H.-D. Cao, B.-L. Chen, and X.-P. Zhu, \emph{Recent developments on {H}amilton's
  {R}icci flow}, Surveys in differential geometry. {V}ol. {XII}. {G}eometric
  flows, vol.~12, Int. Press, Somerville, MA, 2008, pp.~47--112.

\bibitem{caochen}
H.-D. Cao and Q.~Chen, \emph{On locally conformally flat gradient steady
  {R}icci solitons}, Trans. Amer. Math. Soc. \textbf{364} (2012), 2377--2391.

\bibitem{mancat1}
G.~Catino and C.~Mantegazza, \emph{Evolution of the {W}eyl tensor under the
  {R}icci flow}, Ann. Inst. Fourier (2011), 1407--1435.

\bibitem{gahula}
S.~Gallot, D.~Hulin, and J.~Lafontaine, \emph{{R}iemannian geometry},
  Springer--Verlag, 1990.

\bibitem{hamilton1}
R.~S. Hamilton, \emph{Three--manifolds with positive {R}icci curvature}, J.
  Diff. Geom. \textbf{17} (1982), no.~2, 255--306.

\end{thebibliography}

\end{document}